\documentclass[11pt]{article}
\usepackage{amsmath}
\usepackage{amssymb}
\usepackage{amsthm}
\usepackage{subfig}
\usepackage{chngcntr}
\usepackage[centering]{geometry}
\usepackage[dvips]{pict2e}
\usepackage{enumerate}
\usepackage{url}
\usepackage{setspace}

\usepackage{listings}
\usepackage[pdfborder={0 0 0},%
  pdftitle={The minimum rank problem over the finite field of order 2:%
  minimum rank 3},%
  pdfauthor={Wayne Barrett, Jason Grout, and Raphael Loewy},%
  pdfkeywords={Rank 3, Minimum rank, Symmetric matrix, Forbidden subgraph, Field of two elements},
pdfsubject={05C50 (Primary); 05C75, 15A03 (Secondary)}, %
pagebackref=true, %
bookmarksopen=true, bookmarksnumbered=true, bookmarksopenlevel=1, %
pdfpagemode={useOutlines}]{hyperref}

\newcommand{\listindent}{\quad\,}


\newtheoremstyle{proofcases}%
{}{}
{}
{}{\bfseries}{.}
{ }
{\thmname{#1}\thmnumber{ #2}\thmnote{: #3}}

\theoremstyle{plain}
\newtheorem{theorem}{Theorem}
\newtheorem*{theorem*}{Theorem}
\newtheorem{lemma}[theorem]{Lemma}
\newtheorem*{lemma*}{Lemma}
\newtheorem{proposition}[theorem]{Proposition}
\newtheorem{corollary}[theorem]{Corollary}

\theoremstyle{definition}
\newtheorem{definition}[theorem]{Definition}
\newtheorem*{definition*}{Definition}
\newtheorem{observation}[theorem]{Observation}

\newtheorem{example}[theorem]{Example}

\theoremstyle{proofcases}
\newtheorem{case}{Case}
\newtheorem{subcase}{Subcase}[case]

\counterwithin{case}{theorem}

\newcommand{\comp}[1]{\ensuremath{#1^c}}
\newcommand{\setcard}[1]{\ensuremath{\lvert #1 \rvert}}

\newcommand{\st}{\ensuremath{\,\mid\,}}

\newcommand{\I}{\ensuremath{\mathcal{I}}}
\newcommand{\F}{\ensuremath{\mathcal{F}}}
\newcommand{\M}{\ensuremath{\mathcal{MR}}}
\newcommand{\C}{\ensuremath{\mathcal{C}}}
\newcommand{\G}{\ensuremath{\mathcal{G}}}
\newcommand{\0}{\ensuremath{\vec 0}}
\newcommand{\Fld}{\ensuremath{\mathbb{F}}}

\DeclareMathOperator*{\vertexsum}{\oplus}
\newcommand{\vsum}{\vertexsum_v}

\DeclareMathOperator{\mr}{mr}
\DeclareMathOperator{\rank}{rank\,}
\DeclareMathOperator{\col}{col}
\DeclareMathOperator{\row}{row}

\DeclareMathOperator{\wt}{wt}
\DeclareMathOperator{\character}{char}

\newcommand{\bmat}{\begin{bmatrix}}
\newcommand{\emat}{\end{bmatrix}}
\newcommand{\RR}{{\mathbb R}}
\newcommand{\dart}{\mathrm{dart}}
\newcommand{\fh}{\ensuremath{\mathrm{full~house}}}

\newcommand{\keywords}[1]{\noindent\small\textbf{Key words and phrases: #1}}
\newcommand{\subjclass}[1]{\noindent\small\textbf{AMS Subject Classification (2000): #1}}

\begin{document}
\title{The minimum rank problem over the finite field of order 2:
  minimum rank 3}

\author{Wayne Barrett \\ 
Department of Mathematics \\
Brigham Young University \\
Provo, UT 84602, USA\\
\texttt{wayne@math.byu.edu}
\and
Jason Grout\thanks{This work was done at Brigham Young University as part of this author's doctoral dissertation.}\\
Department of Mathematics\\
Iowa State University\\
Ames, IA 50011, USA\\
\texttt{grout@iastate.edu}
\and
Raphael Loewy\thanks{Part of the research of this author was
   done while he was visiting Brigham Young University, whose financial
   support is gratefully acknowledged.}\\
Department of Mathematics\\
Technion-Israel Institute of Technology\\
 Haifa 32000, Israel\\
\texttt{loewy@techunix.technion.ac.il}}

\date{September 2, 2008}

\maketitle

\begin{abstract}
  Our main result is a sharp bound for the number of vertices in
  a minimal forbidden subgraph for  the graphs having minimum
  rank at most 3 over the finite field of order~2.  We also list all
  62 such minimal forbidden subgraphs.  We conclude by exploring how
  some of these results over the finite field of order 2 extend to
  arbitrary fields and demonstrate that at least one third of the 62
  are minimal forbidden subgraphs over an arbitrary field for the
  class of graphs having minimum rank at most 3 in that field.
\end{abstract}

\subjclass{05C50 (Primary); 05C75, 15A03 (Secondary).}

\medskip

\keywords{Rank 3, Minimum rank, Symmetric matrix, Forbidden subgraph,
  Field of two elements.}

\section{Introduction}

Given a field $F$ and a simple undirected graph $G$ on $n$ vertices, let $S(F,G)$ be the
set of symmetric $n\times n$ matrices $A$ with entries in $F$ satisfying
$a_{ij} \neq 0$, $i \neq j$, if and only if $ij$ is an edge in $G$.  There
is no restriction on the diagonal entries of the matrices in $S(F,G)$.
Let
\begin{equation*}
\mr(F,G)=\min\{\rank A \st A \in S(F,G)\}.
\end{equation*}
The problem of finding $\mr(F,G)$ has recently attracted considerable
attention, particularly for the case in which $F = \RR$ (see the
survey paper \cite{fallat-hogben-survey} and the references cited
there or the American Institute of Mathematics workshop website
\cite{aim-workshop-minrank} on this topic).  Relevant papers for us
are \cite{hsieh-minrank, BFH1-minrankpath,
  barrett-vdHL-minrank2-infinite, bento-duarte-tridiag-matrices,
  barrett-vdHL-minrank2-finite, ding-kotlov-minrank-finite}.  In
\cite{barrett-vdHL-minrank2-infinite} and
\cite{barrett-vdHL-minrank2-finite}, the problem of characterizing all
graphs $G$ for which $\mr(F,G) \leq 2$ was addressed.  Complete
characterizations were obtained for all fields and fall into four
cases depending on whether the field is infinite or finite and whether
or not the field characteristic is two. These various classifications
have both striking similarities and distinctive differences.

\begin{figure}[htp]
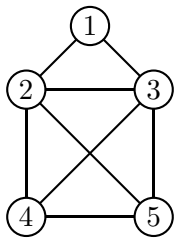

  \centering
  {\rm \input graphs/fullhouse.tex }
  \caption{The full house graph.}
  \label{fig:fullhouse}
\end{figure}

The full house, seen and labeled in Figure~\ref{fig:fullhouse}, is the
only graph on 5 or fewer vertices for which the field affects the
minimum rank.  (This was previously noted in
\cite{barrett-vdHL-minrank2-finite}, in which the graph was identified
as $\overline {P_3 \cup 2K_1}$.)

\begin{proposition}[\cite{barrett-vdHL-minrank2-finite}]\label{prop:no-fh-indep-field}
  Let $G$ be a graph on 5 or fewer vertices and suppose that $G\neq\fh$.
  Then $\mr(F,G)$ is independent of the field $F$.
\end{proposition}

We also include a short proof for the fact that the minimum rank of
the \fh{} graph is field-dependent.

\begin{proposition}
  If $F\neq \Fld_2$ is a field, then $\mr(F,\fh)=2$.  However, $\mr(\Fld_2,\fh)=3$.
\end{proposition}

\begin{proof}
If $F \neq \Fld_2$, there are elements $a, b \neq 0$ in $F$ such that
$a+b \neq 0$.  Then
\begin{equation*}
\begin{bmatrix}
a&a&a&0&0\\
a&a+b&a+b&b&b\\
a&a+b&a+b&b&b\\
0&b&b&b&b\\
0&b&b&b&b
\end{bmatrix} \in S(F,\fh),
\end{equation*}
which shows that $\mr(F,\fh)=2$.  But if
\begin{equation*}
A=  \begin{bmatrix}
d_1&1&1&0&0\\
1&d_2&1&1&1\\
1&1&d_3&1&1\\
0&1&1&d_4&1\\
0&1&1&1&d_5
\end{bmatrix}
\end{equation*}
is any matrix in $S(\Fld_2,\fh)$, then 
\begin{equation*}
\det A[\{1, 2, 5\},\{1, 3,
4\}]=\begin{vmatrix} d_1&1&0\\1&1&1\\0&1&1\end{vmatrix} =1,
\end{equation*}
so $\mr(\Fld_2,\fh)\geq 3$.  Setting all $d_i=1$ verifies that $\mr(\Fld_2,\fh)= 3$.
\end{proof}

In spite of this dependence on the field, it has become clear that
even for calculating the minimum rank over just the real field, 
results obtained over finite fields, and over $\Fld_2$ in
particular, will provide important insights.  This will be explored
more fully in Section~\ref{sec:other-fields}.

The methods of
\cite{barrett-vdHL-minrank2-infinite,barrett-vdHL-minrank2-finite} do
not extend in any straightforward way to the problem of characterizing
graphs with $\mr(F,G) \leq k$ for $k \geq 3$. However, it is possible
to obtain results of this sort for finite fields using other methods
which make explicit use of the finiteness of $F$.  In this
paper we examine the case $F=\Fld_2$.

We recall some notation from graph theory.

\begin{definition} 
  Given a graph $G$, $V(G)$ denotes the set of vertices in $G$ and
  $E(G)$ denotes the set of edges in $G$.  The \emph{order} of a graph is
  $\setcard{G}=\setcard{V(G)}$.  Given two graphs $G$ and $H$, with
  $V(G)$ and $V(H)$ disjoint, the \emph{union} of $G$ and $H$ is $G
  \cup H = (V(G) \cup V(H), E(G) \cup E(H))$.  The \emph{join}, $G
  \lor H$, is the graph obtained from $G \cup H$ by adding edges from
  all vertices of $G$ to all vertices of $H$.  If $S \subset V(G)$,
  $G[S]$ denotes the subgraph of $G$ induced by $S$.  If $H$ is an
  induced subgraph of $G$, $G-H$ denotes the subgraph induced by
  $V(G)\setminus V(H)$.
\end{definition}

\begin{definition} 
  We denote the path on $n$ vertices by $P_n$.  The complete graph on
  $n$ vertices will be denoted by $K_n$.  We abbreviate $K_n \cup
  \cdots \cup K_n$ ($m$ times) to $mK_n$.
\end{definition}

We recall the following observation.
 \begin{observation}[{\cite[Observation 5]{barrett-vdHL-minrank2-infinite}}]
 If $H$ is an induced subgraph of $G$, then for any field
 $F$, $\mr(F,H)\leq\mr(F,G)$.
 \end{observation}

 \begin{example}
   It is well known that $\mr(F,P_{k+2})=k+1$ for any field $F$. Therefore $P_{k+2}$
   cannot be an induced subgraph of any graph $G$ with $\mr(F,G)\leq k$.
\end{example}

\begin{definition}\label{def:F_k+1}
  Let $F$ be any field. The graph $H$ is a \emph{minimal forbidden subgraph} for
  the class of graphs $\G_k (F)=\{G \st \mr(F,G)\leq k\}$ if
  \begin{enumerate}
  \item $\mr(F,H)\geq k+1$ and 
  \item $\mr(F,H-v)\leq k$ for every vertex $v\in V(H)$.
  \end{enumerate}
  Let $\F_{k+1}(F)$ be the set of all minimal forbidden subgraphs for
  $\G_k(F)$.
\end{definition}

 \begin{observation}\label{obs:Gk-iff-Fkplus1-forbidden}
   $G\in \G_k(F)\iff$ no graph in $\F_{k+1}(F)$ is induced in $G$.
 \end{observation}

 Theorem 6 $(a\iff c)$ of \cite{barrett-vdHL-minrank2-infinite} and
 Theorem 16 of \cite{barrett-vdHL-minrank2-finite} can be restated:

 \begin{theorem}[{\cite[Theorem 6]{barrett-vdHL-minrank2-infinite}}]
   \label{thm:F3-R}
   $\F_3(\RR)\ =\{P_4, \ltimes, \dart, P_3\cup K_2, 3K_2,
   K_{3,3,3}\}$.
 \end{theorem}

\begin{theorem}[{\cite[Theorem 16]{barrett-vdHL-minrank2-finite}}]
  \label{thm:F3-F2}
  \begin{equation*}
    \F_3(\Fld_2)=\{P_4, \ltimes, \dart, P_3\cup K_2, 3K_2, \fh,
    P_3\lor P_3\}.
\end{equation*}
 \end{theorem}

 Ding and Kotlov \cite{ding-kotlov-minrank-finite} obtained an
 important result related to $\F_{k+1}(F)$.  They showed that if $F$
 is a finite field, then every graph in $\F_{k+1}(F)$ has at most
 $\left(\frac{|F|^k}{2}+1\right)^2$ vertices, so, in particular,
 $\F_{k+1}(F)$ is finite.  In the special case $k=3$ and $F=\Fld_2$,
 their result implies that each graph in $\F_4(\Fld_2)$ has at most 25
 vertices.  In this paper, we improve their bound for this case to
 show that every graph in $\F_4(\Fld_2)$ has at most 8 vertices.  This
 new bound makes an exhaustive computer search feasible, which gives
 the result that $\setcard{\F_4(\Fld_2)} = 62$ and also shows that the
 new bound is sharp.  Of the 29 graphs in $\F_4(\Fld_2)$ having vertex
 connectivity at most one, we will prove that 21 graphs are in
 $\F_4(F)$ for every field $F$, while none of the remaining 8 graphs
 are in $\F_4(F)$ for any field $F\neq \Fld_2$.

 Our approach relies on the following generalization of $\F_{k+1}(F)$.

\begin{definition}\label{def:F-k+1}
  Given a field $F$ and a graph $H$, let $\F_{k+1}(F,H)$ be the set
  of graphs $G$ containing $H$ as an induced subgraph and
  satisfying
  \begin{enumerate}
  \item $\mr(F,G)\geq k+1$ and\label{item:1}
  \item for some $H$ induced in $G$, $\mr(F,G-v)\leq k$ for every $v\in V(G-H)$.\label{item:2}
  \end{enumerate}
\end{definition}

\begin{example}\label{ex:definition-of-F-k+1}
  Let $F$ be any field, let $G$ be the graph labeled in Figure~\ref{fig:ex-definition},
  \begin{figure}[ht]
    \centering
    {\rm\input graphs/P3cartP2.tex }
    \caption{$G$ in Example~\ref{ex:definition-of-F-k+1}.}
    \label{fig:ex-definition}
  \end{figure}
  let $H=P_4$, and let $k=3$.  Since $P_5$ is induced in $G$,
  $\mr(F,G)\geq 3+1$, so condition~\ref{item:1} is satisfied.

  Six copies of $H=P_4$ are induced in $G$.  For $H=G[\{u,v,w,x\}]$, we
  have $\mr(F,G-y)=4$, so condition~\ref{item:2} is not satisfied
  for this copy of $P_4$.  However, if $H=G[\{u,v,y,z\}]$, both
  $G-w$ and $G-x$ are isomorphic to \input{graphs/P3cartP2minusvertex.tex}, which
  has minimum rank 3 by Theorem 2.3 in \cite{BFH1-minrankpath} (see
  Theorem~\ref{thm:cut-vertex} in this paper).  Therefore
  condition~\ref{item:2} is satisfied for this induced $P_4$,
  so $G\in\F_{k+1}(F,P_4)$.
\end{example}

In the notation of Definition~\ref{def:F-k+1}, $\F_{k+1}(F)=\F_{k+1}(F,\emptyset)$, where
$\emptyset$ is the empty graph.

\begin{theorem}\label{thm:Fk-in-Fkplus1-H}
\begin{equation*}
  \F_{k+1}(F)\subseteq\bigcup_{H\in \F_k(F)} \F_{k+1}(F,H).
\end{equation*}
\end{theorem}
\begin{proof}
  Let $G\in\F_{k+1}(F)$.  Since $\mr(F,G)\geq k+1>k-1$, $G\not\in
  \G_{k-1}(F)$. Therefore some graph $H\in\F_k(F)$ is induced in $G$.
  By definition, $\mr(F,G-v)\leq k$ for every vertex $v$ of $G$, so
  $\mr(F,G-v)\leq k$ for every vertex $v$ of $G-H$.  By definition,
  $G\in\F_{k+1}(F,H)$.
\end{proof}

Combining Theorems~\ref{thm:F3-F2} and \ref{thm:Fk-in-Fkplus1-H}, we
have the following result.
\begin{corollary}\label{cor:F4-in-F3-H}
  \begin{equation*}
    \begin{split}
      \F_4(\Fld_2)&\subseteq\bigcup_{H\in\F_3(\Fld_2)} \F_4(\Fld_2,H) \\
       &\qquad = 
      \F_4(\Fld_2,3K_2) 
      \cup \F_4(\Fld_2,P_3 \lor P_3) 
      \cup \F_4(\Fld_2,\dart)
      \cup \F_4(\Fld_2,\ltimes) \\
      &\qquad\quad \cup \F_4(\Fld_2,P_3\cup K_2) 
      \cup \F_4(\Fld_2,\fh) 
      \cup\F_4(\Fld_2,P_4). 
      \end{split}
  \end{equation*}
\end{corollary}

Sections~\ref{sec:M(F,G)}--\ref{sec:graphs-in-F_4} are devoted to
explicitly determining $\F_4(\Fld_2)$.

\section{Matrices which attain a minimum rank for $\F_3(\Fld_2)$}
\label{sec:M(F,G)}
Given a field $F$ and a graph $G$, it is natural to seek to determine
all matrices in $S(F,G)$ which attain the minimum rank of $G$.
Determining these matrices plays a critical role in determining
$\F_4(\Fld_2)$.

\begin{definition}
  Let $G$ be a graph.  Let $\M(F,G)=\{A\in S(F,G) \st \rank A =
  \mr(F,G)\}$, the set of matrices in $S(F,G)$ that attain the minimum
  rank of $G$.  Call two matrices in $\M(F,G)$ \emph{equivalent} if and only
  if they have the same column space.  Let $\C(F,G)$ be the resulting
  set of equivalence classes.
\end{definition}

Let $G$ be a graph.  In the remainder of this section and in
Sections~\ref{sec:general-theorem}--\ref{sec:p4-lemma}, we will assume
that $F=\Fld_2$ and abbreviate our notation as follows: $S(\Fld_2,G)$
is shortened to $S(G)$, $\mr(\Fld_2,G)$ is shortened to $\mr(G)$,
$\F_{k+1}(\Fld_2)$ is shortened to $\F_{k+1}$, $\F_{k+1}(\Fld_2,G)$ is
shortened to $\F_{k+1}(G)$, $\M(\Fld_2,G)$ is shortened to $\M(G)$,
and $\C(\Fld_2,G)$ is shortened to $\C(G)$.

In the remainder of this section, we determine $\M(G)$ for all of the graphs in $\F_3$ (see Theorem~\ref{thm:F3-F2}).

 \begin{lemma}\label{thm:M-P3} With $P_3$ labeled as \raisebox{-1ex}{{\rm \setlength{\unitlength}{1in}
\begin{picture}(1.2,0.2)(0,-0.1)
\thicklines
\linethickness{0.8pt}
\put(0.1,0){\circle{0.2}}
\put(0.6,0){\circle{0.2}}
\put(1.1,0){\circle{0.2}}
\put(0.2,0){\line(1,0){0.3}}
\put(0.7,0){\line(1,0){0.3}}
\put(0.1,0){\makebox(0,0){1}}
\put(0.6,0){\makebox(0,0){2}}
\put(1.1,0){\makebox(0,0){3}}
\end{picture}
 }},
\begin{equation*}
\M(P_3)=\left\{\bmat
 0 & 1 & 0\\
 1 & 0 & 1\\
 0 & 1 & 0 \emat,
 \bmat
 0 & 1 & 0\\
 1 & 1 & 1\\
 0 & 1 & 0 \emat,
 \bmat
 1 & 1 & 0\\
 1 & 0 & 1\\
 0 & 1 & 1 \emat\right\}.
\end{equation*}
 \end{lemma}

 \begin{proof}
   Since $\mr(P_3)=2$,
\begin{equation*}
A=\bmat
 x & 1 & 0\\
 1 & y & 1\\
 0 & 1 & z \emat
 \in \M(P_3) \iff \det A=xyz+x+z=0\ \text{ in } \Fld_2.
\end{equation*}
If $x\neq z$, then $\det A=1$, so $x=z$. Then $\det A=xy$, so
$A\in \M(P_3)$ if and only if either $x=y=z=0$, $x=z=0$ and $y=1$, or
$x=z=1$ and $y=0$.
\end{proof}

\begin{proposition} \label{thm:M-calc} The sets $\M(G)$ for $G\in\F_3$
  are as follows.
\begin{enumerate}

\item \label{thm:M-3K2} With $3K_2$ labeled as \raisebox{-2.2ex}{{\rm \input graphs/threek2.tex }},
\begin{equation*}
\M(3K_2)=\left\{\bmat
 1 & 1\\
 1 & 1\emat\oplus\bmat
 1 & 1\\
 1 & 1\emat\oplus\bmat
 1 & 1\\
 1 & 1\emat\right\}.
\end{equation*}

\item\label{thm:M-P3VP3} With $P_3\lor P_3$ labeled as
\raisebox{-6ex}{{\rm\input graphs/P3VP3.tex }},
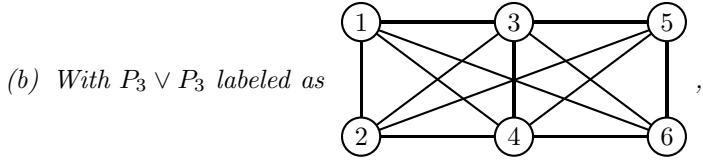
\begin{equation*}
\M(P_3\lor P_3)=\left\{\bmat
 0 & 1 & 1 & 1 & 0 & 1\\
 1 & 0 & 1 & 1 & 1 & 0\\
 1 & 1 & 1 & 1 & 1 & 1\\
 1 & 1 & 1 & 1 & 1 & 1\\
 0 & 1 & 1 & 1 & 0 & 1\\
 1 & 0 & 1 & 1 & 1 & 0\emat\right\}.
\end{equation*}

\item\label{thm:M-dart} With the dart labeled as \raisebox{-7ex}{{\rm \input graphs/dart.tex }},
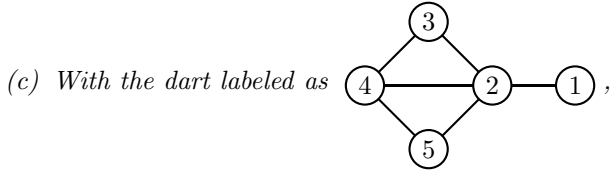
\begin{equation*}
\M(\dart)=\left\{M_1=\bmat
 1 & 1 & 0 & 0 & 0\\
 1 & 0 & 1 & 1 & 1\\
 0 & 1 & 0 & 1 & 0\\
 0 & 1 & 1 & 1 & 1\\
 0 & 1 & 0 & 1 & 0\emat,\, M_2=\bmat
 1 & 1 & 0 & 0 & 0\\
 1 & 1 & 1 & 1 & 1\\
 0 & 1 & 0 & 1 & 0\\
 0 & 1 & 1 & 0 & 1\\
 0 & 1 & 0 & 1 & 0\emat\right\}
\end{equation*} and $\C(\dart)=\{C_1=\{M_1\},\, C_2=\{M_2\}\}$.

\item\label{thm:M-ltimes} With $\ltimes$ labeled as
\raisebox{-7.2ex}{{\rm\input graphs/ltimes.tex }},
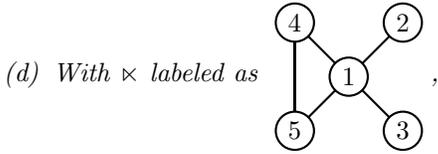
\begin{equation*}
\M(\ltimes)=\left\{M_1=\bmat
 0 & 1 & 1 & 1 & 1\\
 1 & 0 & 0 & 0 & 0\\
 1 & 0 & 0 & 0 & 0\\
 1 & 0 & 0 & 1 & 1\\
 1 & 0 & 0 & 1 & 1\emat,\, M_2=\bmat
 1 & 1 & 1 & 1 & 1\\
 1 & 0 & 0 & 0 & 0\\
 1 & 0 & 0 & 0 & 0\\
 1 & 0 & 0 & 1 & 1\\
 1 & 0 & 0 & 1 & 1\emat,\, M_3=\bmat
 1 & 1 & 1 & 1 & 1\\
 1 & 1 & 0 & 0 & 0\\
 1 & 0 & 1 & 0 & 0\\
 1 & 0 & 0 & 1 & 1\\
 1 & 0 & 0 & 1 & 1\emat\right\}
\end{equation*} and $\C(\ltimes)=\{C_1=\{M_1,M_2\},\, C_2=\{M_3\} \}$.

\item\label{thm:M-P3UK2} With $P_3\cup K_2$ labeled as
  \raisebox{-2.2ex}{{\rm \input graphs/P3UK2.tex }},
\begin{align*} 
\M(P_3\cup K_2)=&\left\{M_1=\bmat
 0 & 1 & 0\\
 1 & 0 & 1\\
 0 & 1 & 0\emat\oplus\bmat
 1 & 1\\
 1 & 1\emat,\,M_2=\bmat
 0 & 1 & 0\\
 1 & 1 & 1\\
 0 & 1 & 0\emat\oplus\bmat
 1 & 1\\
 1 & 1\emat,\right.\\
&\quad\left.M_3=\bmat
 1 & 1 & 0\\
 1 & 0 & 1\\
 0 & 1 & 1\emat\oplus\bmat
 1 & 1\\
 1 & 1\emat\right\}
\end{align*}
and $\C(P_3\cup K_2)=\{C_1=\{M_1,M_2\},\, C_2=\{M_3\}\}$.

\item\label{thm:M-fullhouse} With the full house labeled as in
  Figure~\ref{fig:fullhouse},
 \begin{align*}
 \M(\fh)=&\left\{M_1=\bmat
 1 & 1 & 1 & 0 & 0\\
 1 & 1 & 1 & 1 & 1\\
 1 & 1 & 1 & 1 & 1\\
 0 & 1 & 1 & 1 & 1\\
 0 & 1 & 1 & 1 & 1\emat,\, M_2=\bmat
 0 & 1 & 1 & 0 & 0\\
 1 & 1 & 1 & 1 & 1\\
 1 & 1 & 1 & 1 & 1\\
 0 & 1 & 1 & 1 & 1\\
 0 & 1 & 1 & 1 & 1\emat,\right.\\
 &\quad\left.M_3=\bmat
 0 & 1 & 1 & 0 & 0\\
 1 & 0 & 1 & 1 & 1\\
 1 & 1 & 0 & 1 & 1\\
 0 & 1 & 1 & 1 & 1\\
 0 & 1 & 1 & 1 & 1\emat,\, M_4=\bmat
 1 & 1 & 1 & 0 & 0\\
 1 & 1 & 1 & 1 & 1\\
 1 & 1 & 1 & 1 & 1\\
 0 & 1 & 1 & 0 & 1\\
 0 & 1 & 1 & 1 & 0\emat\right\}
 \end{align*} 
and $\C(\fh)=\{C_1=\{M_1,M_2\},\, C_2=\{M_3\},\, C_3=\{M_4\}\}$.

\item \label{thm:M-P4} With $P_4$ labeled as \raisebox{-2.2ex}{{\rm \input graphs/P4.tex }},
 \begin{align*}
 \M(P_4)=&\left\{M_1=\bmat
 0 & 1 & 0 & 0\\
 1 & 0 & 1 & 0\\
 0 & 1 & 1 & 1\\
 0 & 0 & 1 & 1\emat,\, M_2=\bmat
 0 & 1 & 0 & 0\\
 1 & 1 & 1 & 0\\
 0 & 1 & 1 & 1\\
 0 & 0 & 1 & 1\emat,\right.\\
 &\quad\left.M_3=\bmat
 1 & 1 & 0 & 0\\
 1 & 1 & 1 & 0\\
 0 & 1 & 0 & 1\\
 0 & 0 & 1 & 0\emat,\, M_4=\bmat
 1 & 1 & 0 & 0\\
 1 & 1 & 1 & 0\\
 0 & 1 & 1 & 1\\
 0 & 0 & 1 & 0\emat,\, M_5=\bmat
 1 & 1 & 0 & 0\\
 1 & 0 & 1 & 0\\
 0 & 1 & 0 & 1\\
 0 & 0 & 1 & 1\emat\right\}
\end{align*} and $\C(P_4)=\{C_1=\{M_1,M_2\},\, C_2=\{M_3,M_4\},\,
C_3=\{M_5\}\}$.

\end{enumerate}
 \end{proposition}

\begin{proof}[Proof 1]
  Exhaustively calculate the rank of each matrix in $S(G)$ for each
  $G\in\F_3$.  Appendix~\ref{sec:magma-programs} contains a collection
  of Magma \cite{magma} functions to implement this approach.
\end{proof}

 \begin{proof}[Proof 2] It is known
   \cite{barrett-vdHL-minrank2-infinite,barrett-vdHL-minrank2-finite}
   that each of the graphs in (a) through (g) has minimum rank 3.

Part \ref{thm:M-3K2} follows immediately and \ref{thm:M-P3UK2}
   follows from Lemma~\ref{thm:M-P3}. We prove
   \ref{thm:M-fullhouse} and \ref{thm:M-P4}. The proofs of
   \ref{thm:M-P3VP3}, \ref{thm:M-dart}, and \ref{thm:M-ltimes}
   are similar.
\begin{enumerate}

 \item[\ref{thm:M-fullhouse}] Let
\begin{equation*}
A=\bmat
 v & 1 & 1 & 0 & 0\\
 1 & w & 1 & 1 & 1\\
 1 & 1 & x & 1 & 1\\
 0 & 1 & 1 & y & 1\\
 0 & 1 & 1 & 1 & z\emat
 \in \M(\fh).
\end{equation*}
 \begin{enumerate}[{\listindent}I.]
\item $v=0$. Elementary row and column operations reduce $A$
to
\begin{equation*}
  \bmat
  0 & 1 & 1 & 0 & 0\\
  1 & w & 1 & 0 & 0\\
  1 & 1 & x & 0 & 0\\
  0 & 0 & 0 & y & 1\\
  0 & 0 & 0 & 1 & z\emat .
\end{equation*}
Then we must have $y=z=1$. Since $\begin{vmatrix} 0 & 1 & 1\\1 & w & 1\\1
& 1 &  x\end{vmatrix}=w+x$, we must have $w=x$, so $w=x=1$ or $w=x=0$. This
yields  the matrices $M_2$ and $M_3$ in \ref{thm:M-fullhouse}.

\item $v=1$. Row and column reductions yield
  \begin{equation*}
    B=\bmat
    1 & 0 & 0 & 0 & 0\\
    0 & w+1 & 0 & 1 & 0\\
    0 & 0 & x+1 & 1 & 0\\
    0 & 1 & 1 & y & y+1\\
    0 & 0 & 0 & y+1 & y+z\emat .
  \end{equation*}
  If $w=0$, then $B$ can be further reduced to
  \begin{equation*}
    \bmat
 1 & 0 & 0 & 0 & 0\\
 0 & 1 & 0 & 0 & 0\\
 0 & 0 & x+1 & 1 & 0\\
 0 & 0 & 1 & y+1 & y+1\\
 0 & 0 & 0 & y+1 & y+z\emat,
\end{equation*}
which has rank at least 4, so we must have $w=1$. Then
\begin{equation*}
B=\bmat
 1 & 0 & 0 & 0 & 0\\
 0 & 0 & 0 & 1 & 0\\
 0 & 0 & x+1 & 1 & 0\\
 0 & 1 & 1 & y & y+1\\
 0 & 0 & 0 & y+1 & y+z\emat,
\end{equation*} 
which reduces to
\begin{equation*}
C=\bmat
 1 & 0 & 0 & 0 & 0\\
 0 & 0 & 0 & 1 & 0\\
 0 & 0 & x+1 & 0 & 0\\
 0 & 1 & 0 & 0 & 0\\
 0 & 0 & 0 & 0 & y+z\emat.
\end{equation*}
In order for $\rank C=3$, we require that $x=1$ and $y=z$. This
yields matrices $M_1$ and $M_4$ in \ref{thm:M-fullhouse}.
 \end{enumerate}

  \item[\ref{thm:M-P4}] Let
 \begin{equation*}
A=\bmat
 w & 1 & 0 & 0\\
 1 & x & 1 & 0\\
 0 & 1 & y & 1\\
 0 & 0 & 1 & z\emat
 \in \M(P_4).
\end{equation*}
If $w=0$, by elementary row and column operations the matrix reduces to
\begin{equation*}
  B=\bmat
  0 & 1 & 0 & 0\\
  1 & 0 & 0 & 0\\
  0 & 0 & y & 1\\
  0 & 0 & 1 & z\emat .
\end{equation*}
In order for $\rank B=3$, we must have $y=z=1$, but $x$ can be 0 or
1. This yields matrices $M_1$ and $M_2$ in \ref{thm:M-P4}.  If
$w=1$, one row and column operation gives
\begin{equation*}
  C=\bmat
  1 & 0 & 0 & 0\\
  0 & x+1 & 1 & 0\\
  0 & 1 & y & 1\\
  0 & 0 & 1 & z\emat .
\end{equation*}
In order for $C$ to have rank 3,
\begin{equation*}
  \bmat
  x+1 & 1 & 0\\
  1 & y & 1\\
  0 & 1 & z\emat \text{ must be in } \M(P_3),
\end{equation*}
which by Lemma~\ref{thm:M-P3} gives three possibilities for $x$,
$y$, and $z$, giving matrices $M_3$, $M_4$, and $M_5$ in \ref{thm:M-P4}.

\end{enumerate}

\end{proof}

\section{General theorems}
\label{sec:general-theorem}

Throughout this section, let $G$ be a graph with an induced subgraph
$H$ such that \mbox{$\mr(H)=k$}.

For convenience, in
sections~\ref{sec:general-theorem}--\ref{sec:p4-lemma}, we will
consider $G$ as a complete graph with weighted edges.  The weight of
an edge, $\wt(ij)$, is 1 if $ij$ is an edge in the original graph and 0 if it is not.
The vertices in $G-H$ will also have weights.  Let the vertices of $H$
be labeled $h_1,h_2,\ldots,h_\ell$.  The weight $\wt(v)$ of a vertex
$v\in V(G-H)$ is the vector
$(\wt(vh_1),\wt(vh_2),\ldots,\wt(vh_\ell))^T$ of edge weights between
the vertex $v$ and the vertices of $H$.

\subsection{Definitions}

\begin{definition}
  Let $M\in\M(H)$.  We say the vertex $v$ in $G-H$ is
  \emph{rank-preserving} with respect to $M$ if
\begin{align*}
    \rank\begin{bmatrix} M & \wt(v)\end{bmatrix}  = \rank M.
\end{align*}
If $v$ is rank-preserving with respect to $M$, then $M$ can be augmented by a row
and column to obtain a matrix in $S(G[V(H)\cup\{v\}])$ of rank $k$, so
$\mr(G[V(H)\cup\{v\}])=\mr(H)$.  If $v$ is not rank-preserving with respect to $M$, we
say $v$ is \emph{rank-increasing} with respect to $M$.  We say that a set of
vertices is \emph{rank-preserving} with respect to $M$ if each vertex is rank-preserving
with respect to $M$, and a set is \emph{rank-increasing} with respect to $M$ if some vertex is
rank-increasing with respect to $M$.
\end{definition}

\begin{definition}
  Let $M\in\M(H)$.  We say the edge $uv\in G-H$, $u\neq v$, is
  \emph{rank-preserving} with respect to $M$ if $u$ and $v$ are rank-preserving
  with respect to $M$ and $\wt(uv)$ is the unique number that satisfies the
  equality
  \begin{align*}
    \rank
    \begin{bmatrix}
      M & \wt(u)\\
      \wt(v)^T & \wt(uv)
    \end{bmatrix} = \rank M.
  \end{align*}
  (If $\wt(u)=Mp$ and $\wt(v)=Mq$, then $uv$ is rank-preserving if and
  only if $\wt(uv)=q^TMp$.)  If $uv$ is not rank-preserving with respect to
  $M$, we say $uv$ is \emph{rank-increasing} with respect to $M$.
  Notice that $uv$ is rank-preserving with respect to $M$ if and only if 
  $\mr(G[V(H)\cup\{u,v\}])=\mr(H)$.  We say that a set of edges is \emph{rank-preserving} with respect to
  $M$ if each edge is rank-preserving with respect to $M$ and is \emph{rank-increasing}
  with respect to $M$ if some edge is rank-increasing with respect to $M$.
\end{definition}

We emphasize one part of this definition as:

\begin{observation}\label{thm:edge-vertex-rank-increasing}
  If a vertex $v\in G-H$ is rank-increasing with respect to $M$, then each edge incident
  to $v$ in $G-H$ is also rank-increasing with respect to $M$.
\end{observation}

\begin{definition}
  Let $M\in\M(H)$.  Given an ordered set of vertex weights $v_1,\ldots,v_p \in
  \col(M)$, let $v_i=Ma_i$ and let $A=[a_1\, \cdots \, a_p]$.  Then we say
  that the $p \times p$ matrix $P=A^TMA$ is the \emph{rank-preserving
    table} for the ordered set $v_1,\ldots,v_p$ with respect to $M$.
  Note that the $ij$ entry of $P$ is the edge weight needed to make
  the edge between two vertices with weights $v_i$ and $v_j$ a
  rank-preserving edge with respect to $M$.
\end{definition}

\begin{example}
\label{example:rank-preserving-increasing}
\begin{figure}[htb]
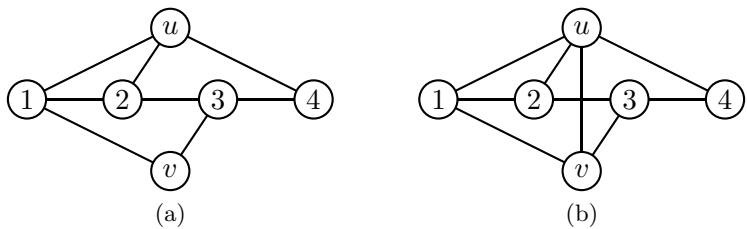

\centering
 \subfloat[][]{\label{fig:example:a}{\rm \input graphs/examplea.tex }}
 \qquad
 \subfloat[][]{\label{fig:example:b}{\rm \input graphs/exampleb.tex }}
  \caption{Graphs in Example \ref{example:rank-preserving-increasing}.}
  \label{fig:example}
\end{figure}
Let $H=P_4$, labeled as in Proposition~\ref{thm:M-calc}, with
corresponding $\M(P_4)$ and $\C(P_4)$.  Let $G$ be a graph containing
vertices $\{1,2,3,4,u,v\}$ such that $H=G[\{1,2,3,4\}]$ and
$G[\{1,2,3,4,u,v\}]$ is one of the graphs in Figure~\ref{fig:example}.
Then $u$ and $v$ have weights $\wt(u)=(1,1,0,1)^T$ and
$\wt(v)=(1,0,1,0)^T$.  The vertex $u$ is rank-preserving with respect
to $M_1$ and $M_2$ since $\wt(u)\in \col(M_1)=\col(M_2)$ and is
rank-increasing with respect to $M_3$, $M_4$, and $M_5$ since
$\wt(u)\not\in \col(M_3)=\col(M_4)$ and $\wt(u)\not\in \col(M_5)$.
Also, $v$ is rank-preserving with respect to $M_1$, $M_2$, and $M_5$
and is rank-increasing with respect to $M_3$ and $M_4$.  The set of
vertices $\{u,v\}$ is rank-preserving with respect to $M_1$ and $M_2$
and is rank-increasing with respect to $M_3$, $M_4$, and $M_5$.

The edge $uv$ is rank-increasing with respect to $M_3$, $M_4$, and
$M_5$ because the set $\{u,v\}$ is rank-increasing with respect to
each of those matrices.  If $G[\{1,2,3,4,u,v\}]$ is the graph in
Figure~\ref{fig:example}\subref{fig:example:a}, then $\wt(uv)=0$ and $uv$ is
rank-preserving with respect to $M_2$ and rank-increasing with respect
to $M_1$.  If $G[\{1,2,3,4,u,v\}]$ is the graph in
Figure~\ref{fig:example}\subref{fig:example:b}, then $\wt(uv)=1$ and $uv$ is
rank-preserving with respect to $M_1$ and rank-increasing with respect
to $M_2$.  Rank-preserving tables with respect to $M_1$ and $M_2$ for
$[\wt(u), \wt(v)]$ are, respectively,
\begin{equation*}
  P_1=\begin{bmatrix}1&1\\1&0\end{bmatrix} \quad\text{and}\quad
  P_2=\begin{bmatrix}0&0\\0&1\end{bmatrix}.
\end{equation*}
\end{example}

Note that $P_1+P_2=J$, the all-ones matrix.  This property will be
important later, so we give it a name now.

\begin{definition}
  Two matrices $A$ and $B$ with entries in $\Fld_2$ are
  \emph{complementary} if $A+B=J$, the all-ones matrix.
\end{definition}

\begin{definition}
  Let $v$ be a vertex in $G-H$ and $V$ be a set of vertices in $G-H$.
  Let
  \begin{equation*}
    \I_v=\{M\in\M(H) \st v \text{ is rank-increasing with respect to } M\}
  \end{equation*} and
  $\I_V=\cup_{v\in V}\I_v$, the set of matrices for which $V$ is
  rank-increasing.  Let 
  \begin{equation*}
    \mathcal{\bar I}_V=\{C \in \mathcal{C} \st V
    \text{ is rank-increasing with respect to every } M\in C\}.    
  \end{equation*}
  
  Let $uv$ be an edge in $G-H$ and $E$ be a set of edges in $G-H$.
  Let
  \begin{equation*}
    \I_{uv}=\{M \in \M(H) \st uv \text{ is rank-increasing with respect to } M\}
  \end{equation*}
  and $\I_E=\cup_{uv\in E}\I_{uv}$, the set of matrices for which $E$ is
  rank-increasing.
\end{definition}

\begin{example}
  We will continue from
  Example~\ref{example:rank-preserving-increasing}. We have
  $\I_u=\{M_3,M_4,M_5\}$ and $\bar \I_u=\{C_2,C_3\}$.  We also have $\I_v=\{M_3,M_4\}$ and
  $\bar \I_v=\{C_2\}$.

  If $\wt(uv)=0$, as is pictured in
  Figure~\ref{fig:example}\subref{fig:example:a}, then
  $\I_{uv}=\{M_1,M_3,M_4,M_5\}$.  If $\wt(uv)=1$, as is pictured in
  Figure~\ref{fig:example}\subref{fig:example:b}, then
  $\I_{uv}=\{M_2,M_3,M_4,M_5\}$.
\end{example}

\begin{observation}\label{obs:wt-in-intersection-col-space}
  Let $V'$ be a set of vertices in $G-H$ such that $\I_{V'}\neq \M(H)$.
  Then for every $v\in V'$,
  \begin{equation*}
    \wt(v)\in \bigcap_{M_i\in\M(H)\setminus \I_{V'}} \col(M_i).
  \end{equation*}

\end{observation}

\subsection{Theorems}

\begin{observation}\label{obs:min-rank-greater}
We have $\mr(G)=k$ if and only if there is
some $M\in\M(H)$ such that every edge and vertex in $G-H$ is rank-preserving
with respect to $M$.  Conversely, $\mr(G)>k$ if and only if there is some set of
edges $E'\subseteq E(G-H)$ and vertices $V'\subseteq V(G-H)$ such that
$\I_{E'}\cup\I_{V'}=\M(H)$.
\end{observation}

Although the inclusion of a vertex set in the second statement is only
necessary in the case when $|G-H|=1$, it will become apparent
that the additional flexibility it introduces enables us to manage
quite a large number of cases in the proof of our main result.

\begin{corollary}\label{cor:proper-subset-neq-MH}
  Assume that $\mr(G)>k$.  If there are sets $E'\subset E(G-H)$ and
  $V'\subset V(G-H)$ such that $\I_{E'}\cup\I_{V'}=\M(H)$ and
  $(\cup_{xy\in E'}\{x,y\})\cup V'\subset V(G-H)$ is a proper subset
  of $V(G-H)$, then $G\not\in\F_{k+1}(H)$.
\end{corollary}

\begin{proof}
  Let $v\in V(G-H)\setminus ((\cup_{xy\in E'}\{x,y\})\cup V')$.  Then $E'\subset E((G-v)-H)$ and $V'\subset V((G-v)-H)$, so $\mr(G-v)>k$ and $G\not\in \F_{k+1}(H)$.
\end{proof}

\begin{proposition} \label{thm:vertex-edge-not-ranki-for-all} Let
  $G\in\F_{k+1}(H)$.  If $\setcard{G-H}\geq 2$, then for every vertex $v$
  in $G-H$, $\I_v\neq \M(H)$.  If $\setcard{G-H}\geq 3$, then for every edge
  $uv$ in $G-H$, $\I_{uv}\neq\M(H)$.
\end{proposition}

\begin{proof}
  Suppose that $G\in\F_{k+1}(H)$.  Suppose there is some
  vertex $v\in G-H$ which is rank-increasing with respect to every
  $M\in\M(H)$.  Let $w$ be a vertex in $G-H$ other than $v$.
  Then $\mr(G-w)>k$, which is a contradiction.

  Similarly, suppose that $G\in\F_{k+1}(H)$.  Suppose there is some edge
  $uv$ in $G-H$ which is rank-increasing with respect to every $M\in
  \M(H)$.  Let $w$ be a vertex in $G-H$ other than $u$ or $v$.  Then
  $\mr(G-w)>k$, which is a contradiction.
\end{proof}

\begin{corollary}\label{thm:vertex-edge-not-ranki-for-all-trivial}
  Let $G\in\F_{k+1}(H)$ and suppose that $\setcard{\M(H)}=1$.  If
  $\setcard{G-H}\geq 2$, then $\I_v=\emptyset$ for every vertex $v$ in $G-H$.  If
  $\setcard{G-H}\geq 3$, then $\I_{uv}=\emptyset$ for every edge $uv$ in $G-H$.
\end{corollary}

\begin{corollary}\label{thm:not-forbidden-if-r-gt-C-and-M-is-one}
  Suppose that $\setcard{\M(H)}=1$.  If
  $G\in\F_{k+1}(H)$, then $\setcard{G-H}\leq 2$.
\end{corollary}

\begin{proof}
  Suppose that $\setcard{\M(H)}=1$ and $\M(H)=\{M\}$.  Then if
  $\setcard{G-H}\geq 3$, $\I_{uv}=\emptyset$ for every edge $uv$ in
  $G-H$ and $\I_v=\emptyset$ for every vertex $v$ in $G-H$.  Since every
  edge and vertex in $G-H$ is rank-preserving with respect to $M$,
  $\mr(G)=\mr(H)=k$ and $G\not\in\F_{k+1}(H)$.
\end{proof}

\begin{corollary}\label{thm:zero-vector}
  Let $G\in\F_{k+1}(H)$.  If $\setcard{G-H}\geq 3$, then 
  $G-H$ contains no vertex $v$ with $\wt(v)=\0$.
\end{corollary}

\begin{proof}
  Let $\setcard{G-H}\geq 3$ and let $v$ be a vertex of $G-H$ with
  $\wt(v)=\0$, the zero vector. 
  Suppose that there is some vertex $w$ of $G-H$ distinct from $v$
  such that the edge $vw$ has nonzero weight.  Then the edge $vw$ is
  rank-increasing for each $M\in \M(H)$, so $\I_{vw}=\M(H)$.  This
  contradicts Proposition~\ref{thm:vertex-edge-not-ranki-for-all}.
  Therefore, $\wt(vw)=0$ for every $w\in V(G-H)$ and $v$ is an
  isolated vertex in $G$.  Therefore $\mr(G)=\mr(G-v)=k$, so
  $G\not\in\F_{k+1}(H)$, a contradiction.
\end{proof}

\begin{lemma}\label{sec:lemma-not-all-increasing}
  Let $G\in\F_{k+1}(H)$.  If $\setcard{G-H}\geq \setcard{\C(H)}+1$,
  then $\I_{V(G-H)}\neq \M(H)$ (i.e., there exists some
  $C\in\mathcal{C}$ for which $V(G-H)$ is rank-preserving with respect to each
  $M\in C$).
\end{lemma}

\begin{proof}
  Suppose that $\I_{V(G-H)}=\M(H)$.  Choose vertices
  $t_1,\ldots,t_{\setcard{\C(H)}}$ from $V(G-H)$ such that $C_i\subseteq\I_{t_i}$
  for $i=1,\ldots,\setcard{\C(H)}$.  Let $T$ be the set containing
  $t_1,\ldots,t_{\setcard{\C(H)}}$.  Then $\setcard{T}\leq \setcard{\C(H)}$
  and $\I_T=\M(H)$.  Let $v\in V(G-H)\setminus T$.  Then $\I_{V(G-H)\setminus\{v\}}=\M(H)$
  and $\mr(G-v)>k$, which is a contradiction.  Thus there is some
  $M\in\M(H)$ and corresponding $C\in\C(H)$ for which $V(G-H)$ is
  rank-preserving.
\end{proof}

By Observation~\ref{obs:min-rank-greater}, $\mr(G)>k$ if and only if
there exist subsets $E'\subseteq E(G-H)$ and $V'\subseteq V(G-H)$
such that $\I_{E'}\cup \I_{V'}=\M(H)$.  We will be interested in
``minimal'' subsets $R\subseteq E(G-H)$ and $T\subseteq V(G-H)$ such
that $\I_R\cup\I_T=\M(H)$ because $R$ and $T$ provide an upper bound for $\setcard{G-H}$, as the following theorem shows.

 \begin{theorem}\label{thm:gen-bounds}
   Assume that $\mr(G)>k$.  Let $R$ be a set of edges in $G-H$ and $T$ be a
   set of vertices in $G-H$ such that $\I_R\cup \I_T = \M(H)$.
   Let $S=\cup_{ij\in R}\{i,j\}$, the set of vertices incident to 
   the edges in $R$.  If $G\in\F_{k+1}(H)$, then $\setcard{G-H}\leq
   \setcard{S}+\setcard{T}\leq 2\setcard{R} + \setcard{T}$.
 \end{theorem}

 \begin{proof}
   We prove the contrapositive.  Suppose that
   $\setcard{G-H}>\setcard{S}+\setcard{T}$ for some $R$, $S$, and $T$
   satisfying the hypotheses.  Let $v\in V(G-H)\setminus \left(S \cup T \right)$ be
   a vertex in $G-H$ that is different from the vertices in $S$ or
   $T$.  Then $\mr(G-v)>k$ and $G\not\in\F_{k+1}(H)$.
 \end{proof}

 The basic idea behind our strategy is to minimize the size of
 $\setcard{S}+\setcard{T}$ to get an upper bound on the number of
 vertices in $G-H$ for which $G\in\F_{k+1}(H)$.

 In our proofs in Sections~\ref{sec:3k2-or-p3Vp3}--\ref{sec:p4-lemma},
 we will examine possible cases for $\I_S$, $\I_R$, and $\I_T$.  The
 following four properties will significantly reduce the number of
 cases we will need to consider.

 Assume that $G$ is a graph such that $\mr(G)>k$.  Let $R\subseteq E(G-H)$ and
 $T\subseteq V(G-H)$.  Let $S=\cup_{ij\in R}\{i,j\}$, the set of
 vertices incident to the edges in $R$.  Then the following properties are
 a direct consequence of the definition of rank-increasing vertices
 and edges.
\begin{enumerate}[{\listindent}P1.]
\item $\I_S$ and $\I_T$ are each the union of equivalence classes in
  $\C(H)$.\label{prop:IS-U-IT-equiv-classes}
\item\label{prop:IS-in-IR} $\I_S\subseteq\I_R$ since if $v\in S$ is
  rank-increasing for a matrix $M\in\M(H)$, then any edge incident to
  $v$ is also rank-increasing for $M$
  (Observation~\ref{thm:edge-vertex-rank-increasing}).
\end{enumerate}

In addition, if $G\in \F_{k+1}(H)$, $\setcard{G-H}\geq
\setcard{\C(H)}+1$, and $\I_R\cup\I_T=\M(H)$, the following properties
are consequences of Lemma~\ref{sec:lemma-not-all-increasing}.
\begin{enumerate}[{\listindent}P1.]
\addtocounter{enumi}{2}
\item $\bar \I_S \cup \bar \I_T \neq \C(H)$.  This implies that $\I_S\neq \M(H)$ and $\I_T\neq \M(H)$.\label{prop:IS-U-IT-neq-MH}
\item There exists a $C\in\C(H)$ such that $C\subseteq\I_R\setminus\I_S$.  This implies that $\I_R\neq \emptyset$.\label{prop:extra-class-in-IR}
\end{enumerate}

Property~\ref{prop:extra-class-in-IR} is a consequence of $\I_R\cup \I_T=\M(H)$ and properties~P\ref{prop:IS-U-IT-equiv-classes} and P\ref{prop:IS-U-IT-neq-MH}.

\begin{definition}
  Assume that $\mr(G)>k$.  Let $\mathcal{A}$ be the set of triples
  $(R,S,T)$ such that
  \begin{enumerate}
  \item $R\subseteq E(G-H)$, $S=\cup_{ij\in R}\{i,j\}$, and
    $T\subseteq V(G-H)$;
  \item $\I_R\cup\I_T=\M(H)$; and
  \item $2\setcard{R}+\setcard{T}$ is minimized.
  \end{enumerate}
  From the triples in $\mathcal{A}$, select those that minimize
  $\setcard{R}$, and from these triples, choose the triples $(R,S,T)$
  that minimize $\setcard{S}$.  We call such an $(R,S,T)$ an
  \emph{optimal triple} for $G$ and $H$.  
\end{definition}

\begin{theorem}\label{thm:minimize-rst}
  Assume that $\mr(G)>k$.  Let $(R,S,T)$ be an optimal triple for $G$ and $H$.
  Then
\begin{enumerate}
\item For every $v\in T$, $\I_v\not\subseteq (\I_{T\setminus \{v\}}\cup\I_{S})$,\label{prop:min-t} and
\item For every $uv\in R$, $\I_{uv}\not\subseteq (\I_{R\setminus \{uv\}}\cup\I_{S}\cup\I_{T})$\label{prop:min-r}.
\end{enumerate}
\end{theorem}

\begin{proof}
  Suppose that $S$ and $T$ do not satisfy \ref{prop:min-t}.  Let
  $v$ be a vertex for which the property does not hold.  Then
  $\I_R\cup\I_{T\setminus\{v\}}=\M(H)$, but
  $2\setcard{R}+\setcard{T\setminus\{v\}}<2\setcard{R}+\setcard{T}$.  This is a
  contradiction since $(R,S,T)\in\mathcal{A}$.

  Suppose that $R$, $S$, and $T$ do not satisfy \ref{prop:min-r}.  Let $uv$ be
  an edge for which the property does not hold.  Let
  $R'=R\setminus\{uv\}$, $S'=\cup_{xy\in R'}\{x,y\}$,
  and $T'=T\cup\{u,v\}$.  Then $\I_{R'}\cup\I_{T'}=\M(H)$ and $2\setcard{R'}+\setcard{T'}\leq 2\setcard{R}+\setcard{T}$, so $(R',S',T')\in\mathcal{A}$.  However, $\setcard{R'}<\setcard{R}$.  This is a
  contradiction since $(R,S,T)$ is an optimal triple.
\end{proof}

The minimality of $|S|$ was not used in the proof of
Theorem~\ref{thm:minimize-rst}, but will be used later.

Let $(R,S,T)$ be an optimal triple for $G$ and $H$.
Theorem~\ref{thm:minimize-rst}\ref{prop:min-t} implies that for every
vertex $v\in T$, there is class of matrices $C\in\C(H)$ such that $v$ is
rank-increasing with respect to every matrix in $C$, while every other vertex in $T$
and every vertex in $S$ is rank-preserving with respect to every matrix in $C$.
Consequently, there are at most $\setcard{\bar \I_T\setminus \bar \I_S}$ vertices
in $T$. 
Theorem~\ref{thm:minimize-rst}\ref{prop:min-r} implies that for every
edge $uv\in R$, there is some matrix $M\in\M(H)$ such that $uv$ is
rank-increasing with respect to $M$, while every other edge in $R$ and
every vertex
in $S\cup T$ is rank-preserving with respect to $M$.  Consequently,
there are at most $\setcard{\I_R\setminus \I_{S\cup T}}=\setcard{\I_R\setminus
(\I_S\cup\I_T)}$ edges in $R$.

\begin{corollary}\label{thm:rt-size}
  Assume that $\mr(G)>k$.  If $(R,S,T)$ is an optimal triple for $G$ and
  $H$, then
  \begin{enumerate}
  \item $\setcard{T}\leq\setcard{\bar \I_T \setminus \bar
      \I_S}=\setcard{\{C\in\C(H) \st
      C\subseteq(\I_T\setminus\I_S)\}}$, and
  \item $\setcard{R}\leq \setcard{\I_R\setminus (\I_S\cup \I_T)}$.
  \end{enumerate}
\end{corollary}

This corollary gives one upper bound for $\setcard{R}$.  There will be
times that we can prove that an edge in $R$ is rank-increasing for one
matrix $M_i\in\M(H)$ if and only if it is also rank-increasing for
another matrix $M_j\in \M(H)$.  In these cases, we can get a smaller
upper bound for $\setcard{R}$.

\begin{corollary}
  Assume that $\mr(G)>k$ and let $(R,S,T)$ be an optimal triple for $G$ and
  $H$.  Then $S\cap T=\emptyset$.
\end{corollary}

\begin{corollary}\label{thm:t-empty}
  Let $G\in\F_{k+1}(H)$ and let $(R,S,T)$ be an optimal triple for $G$ and
  $H$.  If $\I_R=\M(H)$, then $T=\emptyset$ and
  $\setcard{G-H}\leq\setcard{S}\leq 2\setcard{R}$.
\end{corollary}

\begin{proof}
  Since $\I_R=\M(H)$, $\I_R\cup \I_\emptyset=\M(H)$.  Since for any
  $T\subseteq V(G-H)$, 
  $2\setcard{R}+\setcard{\emptyset}\leq 2\setcard{R}+\setcard{T}$, we
  have $T=\emptyset$ by the
  minimality of $2\setcard{R}+\setcard{T}$.  By
  Theorem~\ref{thm:gen-bounds}, $\setcard{G-H}\leq \setcard{S}$.
\end{proof}

The following lemma and corollary give conditions sufficient to reduce the size of
the upper bound for $\setcard{S}$.

\begin{lemma}\label{lem:s-4-to-3-general}
  Assume that $\mr(G)>k$.  Let $(R,S,T)$ be an optimal triple for $G$ and $H$.
  Suppose that
  \begin{enumerate}
  \item $\setcard{R}=2$,
  \item \label{lem:hypothesis-complementary}If $uv$ and $wx$ are any
        two edges between vertices in $S$, then
    either $\I_{wx}\setminus\I_S=\I_{uv}\setminus\I_S$ or
    $\I_{wx}\setminus\I_S=(\I_R\setminus\I_{uv})\setminus\I_S$, and
  \item there are two (not necessarily distinct) vertices $v$ and $w$,
    one incident to each edge of $R$, such that $\I_{\{v,w\}}=\I_S$.
  \end{enumerate}
  Then $\setcard{S}=3$.
\end{lemma}

\begin{proof}
  Since $\setcard{R}=2$, we have $3\leq\setcard{S}\leq4$.  Suppose that
  $\setcard{S}=4$.  Let $R=\{uv,wx\}$ and $S=\{u,v,w,x\}$, where
  $\I_{\{v,w\}}=\I_S$.  Let $A=\I_{uv}\setminus\I_S$ and
  $B=(\I_R\setminus\I_{uv})\setminus\I_S$.  We have $\I_{wx}\setminus\I_S\neq A$ by
  Theorem~\ref{thm:minimize-rst}\ref{prop:min-r}, so
  $\I_{wx}\setminus\I_S=B$ by hypothesis \ref{lem:hypothesis-complementary}.  By hypothesis \ref{lem:hypothesis-complementary}, $\I_{vw}\setminus\I_S=A$ or $\I_{vw}\setminus\I_S=B$.
  \begin{enumerate}[{\listindent}I.]
  \item $\I_{vw}\setminus\I_S=A$.  Let $R'=\{vw,wx\}$ and $S'=\{v,w,x\}$.
  \item $\I_{vw}\setminus\I_S=B$.  Let $R'=\{uv,vw\}$ and $S'=\{u,v,w\}$.
  \end{enumerate}

  Since $\{v,w\}\subseteq S'$, $\I_{S'}=\I_S$.  Also $\I_{R'}=A\cup B
  \cup \I_{S'}=A\cup B\cup \I_{S}=\I_R$.  Therefore, $(R',S',T)$ is a
  triple such that $\I_{R'}\cup\I_T=\M(H)$,
  $2\setcard{R'}+\setcard{T}=2\setcard{R}+\setcard{T}$, and
  $\setcard{R'}=\setcard{R}$, but $\setcard{S'}<\setcard{S}$, which
  contradicts the optimality of $(R,S,T)$.  Thus $\setcard{S}=3$.
\end{proof}

\begin{corollary}\label{cor:s-4-to-3}
  Assume that $\mr(G)>k$.  Let $(R,S,T)$ be an optimal triple for $G$ and $H$.
  Suppose that
  \begin{enumerate}
  \item $\setcard{R}=2$,
  \item If $uv$ and $wx$ are any two edges between vertices in $S$, then
    either $\I_{wx}\setminus\I_S=\I_{uv}\setminus\I_S$ or
    $\I_{wx}\setminus\I_S=(\I_R\setminus\I_{uv})\setminus\I_S$, and
  \item $\setcard{\bar \I_S}\leq 1$.
  \end{enumerate}
  Then $\setcard{S}=3$.
\end{corollary}

\begin{proof}
  Since $\setcard{\bar \I_S}\leq 1$, there is some vertex $y\in S$
  such that $\I_y=\I_S$.  Therefore $\I_{\{y,z\}}= \I_S$ for any
  vertex $z\in S$.  Applying Lemma~\ref{lem:s-4-to-3-general} then gives the result.
\end{proof}

In Sections~\ref{sec:3k2-or-p3Vp3}--\ref{sec:p4-lemma}, we will
determine an upper bound for the number of vertices in graphs in
$\F_4(H)$ for each graph $H$ in
\begin{equation*}
\F_3=\{3K_2, P_3\lor P_3,\dart, \ltimes, P_3\cup K_2, \fh,P_4\}.
\end{equation*}
  We will then apply
Corollary~\ref{cor:F4-in-F3-H} to determine the maximum number of
vertices in a graph in $\F_4$.

\section{$H=3K_2$ or $H=P_3\lor P_3$}
\label{sec:3k2-or-p3Vp3}
By Proposition~\ref{thm:M-calc},
$\setcard{\M(3K_2)}=1$ and $\setcard{\M(P_3\lor P_3)}=1$, so applying 
Corollary~\ref{thm:not-forbidden-if-r-gt-C-and-M-is-one} gives
the following lemma.

\begin{lemma}\label{lem:3k2-p3Vp3}
If $G\in\F_4(3K_2)$ or $G\in\F_4(P_3\lor P_3)$, then $\setcard{G}\leq8$.
\end{lemma}

\section{$H=\text{Dart}$}

\begin{lemma}\label{lem:dart}
  If $G\in\F_4(\dart)$, then $\setcard{G}\leq7$.
\end{lemma}

\begin{proof}
Suppose that $G\in\F_4(\dart)$ and $\setcard{G}\geq 8$ (i.e.,
$\setcard{G-H}\geq 3$).  Then $G-H$ has no vertices with zero weight
by Corollary~\ref{thm:zero-vector}.  Assume that $(R,S,T)$ is an
optimal triple for $G$ and the dart.  Let $\M(\dart)=\{M_1,M_2\}$ and
$\C(\dart)=\{C_1=\{M_1\},C_2=\{M_2\}\}$ be as in
Proposition~\ref{thm:M-calc}\ref{thm:M-dart}.  By property
P\ref{prop:IS-U-IT-equiv-classes}, $\I_S\in\{\emptyset,
C_1,C_2,C_1\cup C_2\}$.  By property P\ref{prop:IS-U-IT-neq-MH},
$\I_S\neq C_1\cup C_2$.  Thus $\I_S\in\{\emptyset, C_1, C_2\}$ and we
have the following cases.

\begin{case}[$\I_S=\emptyset$]\label{sec:restrict-vertex-weights-in-S}
  
By Observation~\ref{obs:wt-in-intersection-col-space}, if $v\in S$,
then 
\begin{equation*}
\wt(v)\in\col(M_1)\cap\col(M_2)=\{\0,\, v_1=(0,1,0,1,0)^T,\,
v_2=(1,0,0,1,0)^T,\, v_3=(1,1,0,0,0)^T\}.  
\end{equation*}
 The rank-preserving tables
for $[v_1,v_2,v_3]$ with respect to $M_1$ and $M_2$ are, respectively,
\begin{align*}
P_1=
\begin{bmatrix}
     0&0&0\\
     0&1&1\\
     0&1&1
\end{bmatrix} \quad \text{ and }\quad
P_2=
\begin{bmatrix}
     0&0&0\\
     0&1&1\\
     0&1&1
\end{bmatrix}.
\end{align*}
Since $P_1=P_2$, an edge in $R$ is rank-preserving for $M_1$ if and
only if it is also rank-preserving for $M_2$.  This combined with
property~P\ref{prop:extra-class-in-IR} implies that
$\I_R=\{M_1,M_2\}$.  Since $M_1\in \I_{uv}$ if and only if $M_2\in
\I_{uv}$ for any edge $uv\in R$,
Theorem~\ref{thm:minimize-rst}\ref{prop:min-r} implies that
$\setcard{R}=1$ and $\setcard{S}=2$.  Since $\I_R=\M(\dart)$,
$T=\emptyset$ and $\setcard{G-H}\leq 2$ by
Corollary~\ref{thm:t-empty}. This contradicts our assumption that
$\setcard{G}\geq8$, so this case cannot occur.
\end{case}
\begin{case}[$\I_S=\{M_1\}$ or $\I_S=\{M_2\}$]

In each of these cases, by property~P\ref{prop:extra-class-in-IR},
$\I_R=\{M_1,M_2\}$.  By Corollary~\ref{thm:rt-size}, $\setcard{R}\leq
1$, so $\setcard{R}=1$.  Again, since $\I_R=\M(\dart)$, $T=\emptyset$
and $\setcard{G-H}\leq 2$ by Corollary~\ref{thm:t-empty}.  This
contradicts our assumption that $\setcard{G}\geq 8$, so neither of
these cases can occur.
\end{case}

Thus $\setcard{G-H}\geq 3$ is impossible, so $\setcard{G-H}\leq 2$ and $\setcard{G}\leq 7$.
\end{proof}

\section{$H=\ltimes$}\label{sec:ltimes}

\begin{lemma}\label{lem:ltimes}
If $G\in\F_4(\ltimes)$, then $\setcard{G}\leq 8$.
\end{lemma}
\begin{proof}
Suppose that $G\in\F_4(\ltimes)$ and $\setcard{G}\geq 8$ (i.e.,
$\setcard{G-H}\geq 3$).  Then $G-H$ has no vertices with zero weight
by Corollary~\ref{thm:zero-vector}.  Assume that $(R,S,T)$ is an
optimal triple for $G$ and $\ltimes$.  Let
$\M(\ltimes)=\{M_1,M_2,M_3\}$ and
$\C(\ltimes)=\{C_1=\{M_1,M_2\},C_2=\{M_3\}\}$ be as in
Proposition~\ref{thm:M-calc}\ref{thm:M-ltimes}.  By property
P\ref{prop:IS-U-IT-equiv-classes}, $\I_S\in\{\emptyset,
C_1,C_2,C_1\cup C_2\}$.  By property P\ref{prop:IS-U-IT-neq-MH},
$\I_S\neq C_1\cup C_2$.  Thus $\I_S\in\{\emptyset, C_1, C_2\}$ and we
have the following cases.

\begin{case}[$\I_S=\emptyset$]

By Observation~\ref{obs:wt-in-intersection-col-space}, if $v\in S$,
then $$\wt(v)\in\bigcap_{i=1}^3\col(M_i)=\{\0,\, v_1=(0,1,1,0,0)^T,\, v_2=(1,0,0,1,1)^T,\, v_3=(1,1,1,1,1)^T\}.$$
The rank-preserving tables for $[v_1,v_2,v_3]$ with respect to $M_1$, $M_2$, and $M_3$ are,
respectively,
\begin{align*}
P_1=\begin{bmatrix}
1&0&1\\
0&1&1\\
1&1&0
\end{bmatrix},\quad
P_2=\begin{bmatrix}
0&0&0\\
0&1&1\\
0&1&1
\end{bmatrix},\quad \text{ and }
P_3=\begin{bmatrix}
0&0&0\\
0&1&1\\
0&1&1
\end{bmatrix}.
\end{align*}
Since $P_2=P_3$, an edge in $R$ is rank-preserving for $M_2$ if and
only if it is also rank-preserving for $M_3$.  Thus, we must either
have both $M_2$ and $M_3$ in $\I_R$ or have neither in the set.  Thus
$\I_R\in\{\{M_1,M_2,M_3\},\, \{M_2,M_3\},\, \{M_1\}\}$.  By property
P\ref{prop:extra-class-in-IR}, $\I_R\neq\{M_1\}$.  Therefore we have
the following cases.

\begin{subcase}[$\I_R=\{M_1,M_2,M_3\}$]
\label{subcase:setup-for-s-3}
Since $\I_R=\M(\ltimes)$, Corollary~\ref{thm:t-empty} implies that
$T=\emptyset$ and $\setcard{G-H}\leq\setcard{S}\leq 2\setcard{R}$.
Since $M_2\in\I_{uv}$ if and only if $M_3\in\I_{uv}$ for each edge
$uv\in R$, $\setcard{R}\leq 2$ by
Theorem~\ref{thm:minimize-rst}\ref{prop:min-r}.  If $\setcard{R}\leq
1$, then $\setcard{G-H}\leq 2$, a contradiction.  If $\setcard{R}=2$,
then Theorem~\ref{thm:minimize-rst}\ref{prop:min-r} implies that $R$
consists of an edge $uv$ such that $\I_{uv}=\{M_1\}$ and another edge
$wx$ such that $\I_{wx}=\{M_2,M_3\}$.  Since the second row and column
of $P_1$, $P_2$, and $P_3$ are identical, we see that if any vertex in
$S$, say $u$, has weight $v_2$, then the edge in $R$ incident to the
vertex must have either $\I_{uv}=\I_R$ or $\I_{uv}=\emptyset$.
Neither of these cases occur, so $u$, $v$, $w$, and $x$ each must have
weight $v_1$ or $v_3$.  Note that since the principal submatrices
$P_1[1,3]$ and $P_2[1,3]=P_3[1,3]$ are complementary, any edge between
vertices with weights $v_1$ or $v_3$ must be either rank-increasing
for $M_1$ and rank-preserving for $M_2$ and $M_3$, or rank-increasing
for $M_2$ and $M_3$ and rank-preserving for $M_1$.  This fact combined
with the facts that $\setcard{R}=2$ and $\setcard{\bar \I_S}=0$ allow
us to apply Corollary~\ref{cor:s-4-to-3} to conclude that
$\setcard{S}=3$, $\setcard{G-H}\leq3$, and $\setcard{G}\leq8$.
\end{subcase}
\begin{subcase}[$\I_R=\{M_2,M_3\}$]
Since $\I_R\cup\I_T=\M(H)$,
properties~P\ref{prop:IS-U-IT-equiv-classes} and
P\ref{prop:IS-U-IT-neq-MH} imply that $\I_T=\{M_1,M_2\}$.  By
Corollary~\ref{thm:rt-size}, $\setcard{R}\leq 1$ and $\setcard{T}\leq
1$. By Theorem~\ref{thm:gen-bounds}, $\setcard{G-H}\leq 3$, so
$\setcard{G}\leq 8$.
\end{subcase}
\end{case}
\begin{case}[$\I_S=\{M_1,M_2\}$]
By property~P\ref{prop:extra-class-in-IR}, $\I_R=\{M_1,M_2,M_3\}$, so
$T=\emptyset$ and $\setcard{G-H}\leq 2\setcard{R}$ by
Corollary~\ref{thm:t-empty}.  By Corollary~\ref{thm:rt-size},
$\setcard{R}\leq 1$, so $\setcard{G-H}\leq 2$.  This contradicts the
assumption that $\setcard{G-H}\geq 3$, so this case cannot occur.
\end{case}
\begin{case}[$\I_S=\{M_3\}$]

By Observation~\ref{obs:wt-in-intersection-col-space}, if $v\in S$,
then 
\begin{equation*}
  \begin{split}
 \wt(v)\in\col(M_1)=\col(M_2)
 =\{\0,\, v_1=(0,0,0,1,1)^T,\, v_2=(0,1,1,0,0)^T,\, v_3=(0,1,1,1,1)^T,\,\\
 v_4=(1,0,0,0,0)^T,\, v_5=(1,0,0,1,1)^T,\, v_6=(1,1,1,0,0)^T,\, v_7=(1,1,1,1,1)^T\}. 
  \end{split}
\end{equation*}
The rank-preserving tables for $[v_1,v_2,v_3,v_4,v_5,v_6,v_7]$ with
respect to $M_1$ and $M_2$ are, respectively,
\begin{align*}
P_1=\begin{bmatrix}
1&1&0&0&1&1&0\\
1&1&0&1&0&0&1\\
0&0&0&1&1&1&1\\
0&1&1&0&0&1&1\\
1&0&1&0&1&0&1\\
1&0&1&1&0&1&0\\
0&1&1&1&1&0&0
\end{bmatrix}\quad \text{ and } \quad
P_2=\begin{bmatrix}
1&1&0&0&1&1&0\\
1&0&1&1&0&1&0\\
0&1&1&1&1&0&0\\
0&1&1&0&0&1&1\\
1&0&1&0&1&0&1\\
1&1&0&1&0&0&1\\
0&0&0&1&1&1&1
\end{bmatrix}.
\end{align*}
By property~P\ref{prop:extra-class-in-IR}, $\I_R=\{M_1,M_2,M_3\}$, so
$T=\emptyset$ and $\setcard{G-H}\leq \setcard{S}$ by
Corollary~\ref{thm:t-empty}.  By Corollary~\ref{thm:rt-size},
$\setcard{R}\leq 2$.  If $\setcard{R}=1$, then $\setcard{G-H}\leq 2$,
which is a contradiction.  If $\setcard{R}=2$, then
Theorem~\ref{thm:minimize-rst}\ref{prop:min-r} implies that $R$
consists of an edge $uv$ such that $\I_{uv}\setminus \I_S=\{M_1\}$ and
another edge $wx$ such that $\I_{wx}\setminus\I_S=\{M_2\}$.  Since the
first, fourth, and fifth rows and columns of $P_1$ and $P_2$ are
identical, we see that if any vertex, say $u$, has weight $v_1$,
$v_4$, or $v_5$, then the edge in $R$ incident to the vertex must have
either $\I_{uv}\setminus\I_S=\I_R\setminus\I_S$ or $\I_{uv}\setminus
\I_S=\emptyset$.  Neither of these cases occur, so $u$, $v$, $w$, and
$x$ each must have weight $v_2$, $v_3$, $v_6$, or $v_7$.  As in
Subcase~\ref{subcase:setup-for-s-3}, since $P_1[2,3,6,7]$ and $P_2[2,3,6,7]$ are
complementary, $\setcard{R}=2$, and $\setcard{\bar \I_S}=1$, we can
apply Corollary~\ref{cor:s-4-to-3} to conclude that $\setcard{S}=3$,
$\setcard{G-H}\leq 3$, and $\setcard{G}\leq 8$.\qedhere
\end{case}
\end{proof}

\section{$H=P_3\cup K_2$}

\begin{lemma}\label{lem:p3Uk2}
  If $G\in\F_4(P_3\cup K_2)$, then $\setcard{G}\leq 8$.
\end{lemma}
\begin{proof}
Suppose that $G\in \F_4(P_3\cup K_2)$ and $\setcard{G}\geq 8$ (i.e., $\setcard{G-H}\geq 3$).
  Then $G-H$ has no vertices with zero weight by
Corollary~\ref{thm:zero-vector}.  Assume that $(R,S,T)$ is an optimal
triple for $G$ and $P_3\cup K_2$.  Let $\M(P_3\cup K_2)=\{M_1,M_2,M_3\}$ and $\C(P_3\cup
K_2)=\{C_1=\{M_1,M_2\},C_2=\{M_3\}\}$ be as in
Proposition~\ref{thm:M-calc}\ref{thm:M-P3UK2}. By
properties~P\ref{prop:IS-U-IT-equiv-classes} and
P\ref{prop:IS-U-IT-neq-MH}, $\I_S\in\{\emptyset, C_1,C_2\}$, so we have
the following cases.

\begin{case}[$\I_S=\emptyset$]
By Observation~\ref{obs:wt-in-intersection-col-space}, if $v\in S$,
then $$\wt(v)\in\bigcap_{i=1}^3\col(M_i)=\{\0,\, v_1=(0,0,0,1,1)^T,\, v_2=(1,0,1,0,0)^T,\, v_3=(1,0,1,1,1)^T\}.$$  The rank-preserving
tables for $[v_1,v_2,v_3]$ with respect to $M_1$, $M_2$, and $M_3$ are,
respectively,
\begin{align*}
P_1=\begin{bmatrix}
1&0&1\\
0&0&0\\
1&0&1
\end{bmatrix},\quad
P_2=\begin{bmatrix}
1&0&1\\
0&1&1\\
1&1&0
\end{bmatrix},\quad \text{ and }
P_3=\begin{bmatrix}
1&0&1\\
0&0&0\\
1&0&1
\end{bmatrix}.
\end{align*}
Since $P_1=P_3$, an edge in $R$ is rank-preserving for $M_1$ if and
only if it is also rank-preserving for $M_3$.  Thus
$\I_R\in\{\{M_1,M_2,M_3\},\, \{M_1,M_3\},\, \{M_2\}\}$.  By property~P\ref{prop:extra-class-in-IR}, $\I_R\neq \{M_2\}$.  Therefore we have
the following cases.

\begin{subcase}[$\I_R=\{M_1,M_2,M_3\}$]
Since $\I_R=\M(P_3\cup K_2)$, Corollary~\ref{thm:t-empty} implies that
$T=\emptyset$ and $\setcard{G-H}\leq \setcard{S}\leq 2\setcard{R}$.  We
reason as in Subcase~\ref{subcase:setup-for-s-3} in Section~\ref{sec:ltimes}. Since $M_1\in\I_{uv}$ if
and only if $M_3\in\I_{uv}$ for each edge $uv\in R$, $\setcard{R}\leq
2$ by Theorem~\ref{thm:minimize-rst}\ref{prop:min-r}.  If
$\setcard{R}=1$, then $\setcard{G-H}\leq 2$, a contradiction.  If
$\setcard{R}=2$, then Theorem~\ref{thm:minimize-rst}\ref{prop:min-r}
implies that $R$ consists of an edge $uv$ such that
$\I_{uv}=\{M_1,M_3\}$ and another edge $wx$ such that
$\I_{wx}=\{M_2\}$.  Since the first row and column of $P_1$, $P_2$,
and $P_3$ are identical, we see that if any vertex, say $u$, has
weight $v_1$, then the edge in $R$ incident to the vertex must have
either $\I_{uv}=\I_R$ or $\I_{uv}=\emptyset$.  Neither of these cases
occur, so $u$, $v$, $w$, and $x$ each must have weight $v_2$ or $v_3$.
As in Subcase~\ref{subcase:setup-for-s-3} in Section~\ref{sec:ltimes}, since $P_1[2,3]$ and $P_2[2,3]$ are
complementary, $\setcard{R}=2$, and $\setcard{\bar \I_S}=0$, we can
apply Corollary~\ref{cor:s-4-to-3} to conclude that $\setcard{S}=3$,
$\setcard{G-H}\leq 3$, and $\setcard{G}\leq 8$.

\end{subcase}
\begin{subcase}[$\I_R=\{M_1,M_3\}$]
Since $\I_R\cup\I_T=\M(H)$,
properties~P\ref{prop:IS-U-IT-equiv-classes} and
P\ref{prop:IS-U-IT-neq-MH} imply that $\I_T=\{M_1,M_2\}$.  By
Corollary~\ref{thm:rt-size}, $\setcard{R}\leq 1$ and $\setcard{T}\leq
1$.  By Theorem~\ref{thm:gen-bounds}, $\setcard{G-H}\leq 3$, so
$\setcard{G}\leq 8$.
\end{subcase}
\end{case}
\begin{case}[$\I_S=\{M_1,M_2\}$]
By property~P\ref{prop:extra-class-in-IR}, $\I_R=\{M_1,M_2,M_3\}$, so
$T=\emptyset$ and $\setcard{G-H}\leq 2\setcard{R}$ by
Corollary~\ref{thm:t-empty}.  By Corollary~\ref{thm:rt-size},
$\setcard{R}\leq 1$, so $\setcard{G-H}\leq 2$ and $\setcard{G}\leq 7$.
This contradicts the assumption that $\setcard{G}\geq 8$, so this case
cannot occur.
\end{case}
\begin{case}[$\I_S=\{M_3\}$]
By Observation~\ref{obs:wt-in-intersection-col-space}, if $v\in S$,
then 
\begin{equation*}
  \begin{split}
\wt(v)\in\col(M_1)=\col(M_2)=\{\0,\, v_1=(0,0,0,1,1)^T,\,
v_2=(0,1,0,0,0)^T,\, v_3=(0,1,0,1,1)^T,\,\\ v_4=(1,0,1,0,0)^T,\,
v_5=(1,0,1,1,1)^T,\, v_6=(1,1,1,0,0)^T,\, v_7=(1,1,1,1,1)^T\}.
  \end{split}
\end{equation*}
The
rank-preserving tables for $[v_1,v_2,v_3,v_4,v_5,v_6,v_7]$ with respect to $M_1$ and
$M_2$ are, respectively,
\begin{align*}
P_1=\begin{bmatrix}
     1&0&1&0&1&0&1\\
     0&0&0&1&1&1&1\\
     1&0&1&1&0&1&0\\
     0&1&1&0&0&1&1\\
     1&1&0&0&1&1&0\\
     0&1&1&1&1&0&0\\
     1&1&0&1&0&0&1
\end{bmatrix}\quad \text{ and } \quad
P_2=\begin{bmatrix}
     1&0&1&0&1&0&1\\
     0&0&0&1&1&1&1\\
     1&0&1&1&0&1&0\\
     0&1&1&1&1&0&0\\
     1&1&0&1&0&0&1\\
     0&1&1&0&0&1&1\\
     1&1&0&0&1&1&0
\end{bmatrix}.
\end{align*}
By property~P\ref{prop:extra-class-in-IR}, $\I_R=\{M_1,M_2,M_3\}$, so
$T=\emptyset$ and $\setcard{G-H}\leq\setcard{S}$ by
Corollary~\ref{thm:t-empty}.  By Corollary~\ref{thm:rt-size},
$\setcard{R}\leq 2$.  If $\setcard{R}=1$, then $\setcard{G-H}\leq 2$,
which is a contradiction.  If $\setcard{R}=2$, then
Theorem~\ref{thm:minimize-rst}\ref{prop:min-r} implies that $R$
consists of an edge $uv$ such that $\I_{uv}\setminus \I_S=\{M_1\}$ and
another edge $wx$ such that $\I_{wx}\setminus\I_S=\{M_2\}$.  Since the
first three rows and columns of $P_1$ and $P_2$ are identical, we see
that if any vertex, say $u$, has weight $v_1$, $v_2$, or $v_3$, then
the edge in $R$ incident to the vertex must have either
$\I_{uv}\setminus\I_S=\I_R\setminus\I_S$ or $\I_{uv}\setminus
\I_S=\emptyset$.  Neither of these cases occur, so $u$, $v$, $w$, and
$x$ each must have weight $v_4$, $v_5$, $v_6$, or $v_7$.  As in
Subcase~\ref{subcase:setup-for-s-3} in Section~\ref{sec:ltimes}, since $P_1[4,5,6,7]$ and $P_2[4,5,6,7]$ are
complementary, $\setcard{R}=2$, and $\setcard{\bar \I_S}=1$, we can
apply Corollary~\ref{cor:s-4-to-3} to conclude that $\setcard{S}=3$,
$\setcard{G-H}\leq 3$, and $\setcard{G}\leq 8$. \qedhere
\end{case}
\end{proof}

\section{$H=\fh$}

\begin{lemma}\label{lem:fullhouse}
  If $G\in\F_4(\fh)$, then $\setcard{G}\leq 8$.
\end{lemma}
\begin{proof}
Suppose that $G\in\F_4(\fh)$ and $\setcard{G}\geq 9$ (i.e.,
$\setcard{G-H}\geq 4$).  Then $G-H$ has no vertices with zero weight
by Corollary~\ref{thm:zero-vector}.  Assume that $(R,S,T)$ is an
optimal triple for $G$ and the $\fh$.  Let
$\M(\fh)=\{M_1,M_2,M_3,M_4\}$ and
$\C(\fh)=\{C_1=\{M_1,M_2\},C_2=\{M_3\},C_3=\{M_4\}\}$ be as in
Proposition~\ref{thm:M-calc}\ref{thm:M-fullhouse}.  By
properties~P\ref{prop:IS-U-IT-equiv-classes} and
P\ref{prop:IS-U-IT-neq-MH}, $\I_S\in\{\emptyset, C_1,C_2,C_3,C_1\cup
C_2, C_1\cup C_3, C_2\cup C_3\}$, so we
have the following cases.

\begin{case}[$\I_S=\emptyset$]
By Observation~\ref{obs:wt-in-intersection-col-space}, if $v\in S$,
then $$\wt(v)\in\bigcap_{i=1}^4\col(M_i)=\{\0,\, v_1=(0,0,0,1,1)^T\}.$$ The
rank-preserving tables for $[v_1]$ with respect to $M_1$, $M_2$, $M_3$,
and $M_4$ are, respectively,
\begin{align*}
P_1=
\begin{bmatrix}
  0
\end{bmatrix}, \quad
P_2=
\begin{bmatrix}
  1
\end{bmatrix}, \quad
P_3=
\begin{bmatrix}
  1
\end{bmatrix}, \quad \text{ and }
P_4=
\begin{bmatrix}
  0
\end{bmatrix}.
\end{align*}
Since $P_1=P_4$ and $P_2=P_3$ are complementary, an edge $uv$ in $R$
has either $\I_{uv}=\{M_1,M_4\}$ or $\I_{uv}=\{M_2,M_3\}$.  Thus
$\I_R\in\{\{M_1,M_2,M_3,M_4\},\, \{M_1,M_4\},\, \{M_2,M_3\}\}$.
Therefore we have the following cases.
\begin{subcase}[$\I_R=\{M_1,M_2,M_3,M_4\}$]
Since $\I_R=\M(\fh)$, Corollary~\ref{thm:t-empty} implies that
$T=\emptyset$ and $\setcard{G-H}\leq \setcard{S}\leq 2\setcard{R}$.
We reason as in Subcase~\ref{subcase:setup-for-s-3} in Section~\ref{sec:ltimes}. Since $M_1\in\I_{uv}$
if and only if $M_4\in\I_{uv}$ and $M_2\in\I_{uv}$ if and only if
$M_3\in\I_{uv}$ for each edge $uv\in R$, $\setcard{R}\leq 2$ by
Theorem~\ref{thm:minimize-rst}\ref{prop:min-r}.  If $\setcard{R}=1$,
then $\setcard{G-H}\leq 2$, a contradiction.  If $\setcard{R}=2$, then
Theorem~\ref{thm:minimize-rst}\ref{prop:min-r} implies that $R$
consists of an edge $uv$ such that $\I_{uv}=\{M_1,M_4\}$ and another
edge $wx$ such that $\I_{wx}=\{M_2,M_3\}$.  As in
Subcase~\ref{subcase:setup-for-s-3} in Section~\ref{sec:ltimes}, since $P_1=P_4$ and $P_2=P_3$ are
complementary, $\setcard{R}=2$, and $\setcard{\bar \I_S}=0$, we can
apply Corollary~\ref{cor:s-4-to-3} to conclude that $\setcard{S}=3$,
 $\setcard{G-H}\leq 3$, and $\setcard{G}\leq 8$.  This contradicts
the assumption that $\setcard{G}\geq 9$, so this case does not occur.

\end{subcase}
\begin{subcase}[$\I_R=\{M_1,M_4\}$]
\label{subcase:fullhouse-contradiction}
Since $\I_R\cup\I_T=\M(H)$,
properties~P\ref{prop:IS-U-IT-equiv-classes} and
P\ref{prop:IS-U-IT-neq-MH} imply that $\I_T=\{M_1,M_2,M_3\}$.  By
Corollary~\ref{thm:rt-size}, $\setcard{R}\leq 1$ and $\setcard{T}\leq
2$.  If $\setcard{T}=2$, then by
Theorem~\ref{thm:minimize-rst}\ref{prop:min-t}, $T$ consists of a
vertex $v$ such that $\I_v=\{M_1,M_2\}$ and another vertex $w$ such
that $\I_w=\{M_3\}$.  Since $v$ is rank-preserving with respect to
$M_3$ and $M_4$, $\wt(v)\in \col(M_3)\cap \col(M_4)=\{\0,\,
(0,0,0,1,1)^T\}$, so $\wt(v)=(0,0,0,1,1)^T$.  But then
$\I_v=\emptyset$, a contradiction, so this case does not occur.

\end{subcase}
\begin{subcase}[$\I_R=\{M_2,M_3\}$]
Since $\I_R\cup\I_T=\M(H)$, properties~P\ref{prop:IS-U-IT-equiv-classes} and P\ref{prop:IS-U-IT-neq-MH} imply that $\I_T=\{M_1,M_2,M_4\}$.
By Corollary~\ref{thm:rt-size}, $\setcard{R}\leq 1$ and
$\setcard{T}\leq 2$.  Again, if $\setcard{T}=2$, then by
Theorem~\ref{thm:minimize-rst}\ref{prop:min-t}, $T$ consists of a
vertex $v$ such that $\I_v=\{M_1,M_2\}$ and another vertex $w$ such
that $\I_w=\{M_4\}$.  Proceeding as in
Subcase~\ref{subcase:fullhouse-contradiction}, $\wt(v)=(0,0,0,1,1)^T$ and
$\I_v=\emptyset$, a contradiction, so this case does not occur.
\end{subcase}
\end{case}
\begin{case}[$\I_S=\{M_1,M_2\}$]
By Observation~\ref{obs:wt-in-intersection-col-space}, if $v\in S$,
then $\wt(v)\in\col(M_3)\cap\col(M_4)=\{\0,\, v_1=(0,0,0,1,1)^T\}$.  But
then $\I_S=\emptyset$, a contradiction, so this case does not occur.

\end{case}
\begin{case}[$\I_S=\{M_3\}$]
By Observation~\ref{obs:wt-in-intersection-col-space}, if $v\in S$,
then
\begin{equation*}
  \begin{split}
    \wt(v)\in\col(M_1)\cap\col(M_2)\cap\col(M_4)=\{\0,\,
    v_1=(0,0,0,1,1)^T,\, v_2=(1,1,1,0,0)^T,\,\\ v_3=(1,1,1,1,1)^T\}.
  \end{split}
\end{equation*}
The rank-preserving tables for $[v_1,v_2,v_3]$ with respect to $M_1$,
$M_2$, and $M_4$ are, respectively,
\begin{align*}
P_1=\begin{bmatrix}
0&0&0\\
0&1&1\\
0&1&1
\end{bmatrix},\quad
P_2=\begin{bmatrix}
1&1&0\\
1&0&1\\
0&1&1
\end{bmatrix},\quad \text{ and }
P_4=\begin{bmatrix}
0&0&0\\
0&1&1\\
0&1&1
\end{bmatrix}.
\end{align*}
Since $P_1=P_4$, an edge in $R$ is rank-preserving for $M_1$ if and
only if it is also rank-preserving for $M_4$.  Thus
$\I_R\in\{\{M_1,M_2,M_3,M_4\},\, \{M_1,M_3,M_4\},\, \{M_2,M_3\}\}$.  By property~P\ref{prop:extra-class-in-IR}, $\I_R\neq \{M_2,M_3\}$.  Therefore we have
the following cases.

\begin{subcase}[$\I_R=\{M_1,M_2,M_3,M_4\}$]
Since $\I_R=\M(\fh)$, Corollary~\ref{thm:t-empty} implies that
$T=\emptyset$ and $\setcard{G-H}\leq \setcard{S}$.  We reason as in 
Subcase~\ref{subcase:setup-for-s-3} in Section~\ref{sec:ltimes}. Since $M_1\in\I_{uv}$ if and only if
$M_4\in\I_{uv}$ for each edge $uv\in R$, $\setcard{R}\leq 2$ by
Theorem~\ref{thm:minimize-rst}\ref{prop:min-r}.  If $\setcard{R}=1$,
then $\setcard{G-H}\leq 2$, which is a contradiction.  If
$\setcard{R}=2$, then Theorem~\ref{thm:minimize-rst}\ref{prop:min-r}
implies that $R$ consists of an edge $uv$ such that
$\I_{uv}\setminus\I_S=\{M_1,M_4\}$ and another edge $wx$ such that
$\I_{wx}\setminus\I_S=\{M_2\}$.  Since the third row and column of
$P_1$, $P_2$, and $P_4$ are identical, we see that if any vertex, say
$u$, has weight $v_3$, then the edge in $R$ incident to the vertex
must have either $\I_{uv}\setminus\I_S=\I_R\setminus\I_S$ or
$\I_{uv}\setminus\I_S=\emptyset$.  Neither of these cases occur, so
$u$, $v$, $w$, and $x$ each must have weight $v_1$ or $v_2$.  As in
Subcase~\ref{subcase:setup-for-s-3} in Section~\ref{sec:ltimes}, since $P_1[1,2]=P_4[1,2]$ and $P_2[1,2]$ are
complementary, $\setcard{R}=2$, and $\setcard{\bar \I_S}=1$, we can
apply Corollary~\ref{cor:s-4-to-3} to conclude that $\setcard{S}=3$,
$\setcard{G-H}\leq 3$, and $\setcard{G}\leq 8$.   This contradicts the
assumption that $\setcard{G}\geq 9$, so this case does not occur.

\end{subcase}
\begin{subcase}[$\I_R=\{M_1,M_3,M_4\}$]
Since $\I_R\cup\I_T=\M(H)$,
properties~P\ref{prop:IS-U-IT-equiv-classes} and
P\ref{prop:IS-U-IT-neq-MH} imply that $\I_T=\{M_1,M_2\}$ or
$\I_T=\{M_1,M_2,M_3\}$.  In each of these cases, $\setcard{R}\leq1$
and $\setcard{T}\leq1$ by
Theorem~\ref{thm:minimize-rst}\ref{prop:min-r} and
Corollary~\ref{thm:rt-size}, implying that $\setcard{G-H}\leq3$ and
$\setcard{G}\leq 8$.  This contradicts the assumption that
$\setcard{G}\geq 9$, so these cases do not occur.
\end{subcase}
\end{case}
\begin{case}[$\I_S=\{M_4\}$]
By Observation~\ref{obs:wt-in-intersection-col-space}, if $v\in S$,
then
\begin{equation*}
  \begin{split}
    \wt(v)\in\col(M_1)\cap\col(M_2)\cap\col(M_3)=\{\0,\,
    v_1=(0,0,0,1,1)^T,\, v_2=(0,1,1,0,0)^T,\,\\ v_3=(0,1,1,1,1)^T\}.
  \end{split}
\end{equation*}
The rank-preserving tables for $[v_1,v_2,v_3]$ with respect to $M_1$,
$M_2$, and $M_3$ are, respectively,
\begin{align*}
P_1=\begin{bmatrix}
0&1&1\\
1&1&0\\
1&0&1
\end{bmatrix},\quad
P_2=\begin{bmatrix}
1&0&1\\
0&0&0\\
1&0&1
\end{bmatrix},\quad \text{ and }
P_3=\begin{bmatrix}
1&0&1\\
0&0&0\\
1&0&1
\end{bmatrix}.
\end{align*}
Since $P_2=P_3$, an edge in $R$ is rank-preserving for $M_2$ if and
only if it is also rank-preserving for $M_3$.  By
properties~P\ref{prop:IS-in-IR} and P\ref{prop:extra-class-in-IR},
$\I_R\in\{\{M_1,M_2,M_3,M_4\},\, \{M_2,M_3,M_4\}\}$, so we have the
following cases.

\begin{subcase}[$\I_R=\{M_1,M_2,M_3,M_4\}$]
Since $\I_R=\M(\fh)$, Corollary~\ref{thm:t-empty} implies that
$T=\emptyset$ and $\setcard{G-H}\leq \setcard{S}$.
We again reason as in Subcase~\ref{subcase:setup-for-s-3} in Section~\ref{sec:ltimes}. Since $M_2\in\I_{uv}$
if and only if $M_3\in\I_{uv}$ for each edge $uv\in R$,
$\setcard{R}\leq 2$ by
Theorem~\ref{thm:minimize-rst}\ref{prop:min-r}.  If $\setcard{R}=1$,
then $\setcard{G-H}\leq 2$, which is a contradiction.  If
$\setcard{R}=2$, then Theorem~\ref{thm:minimize-rst}\ref{prop:min-r}
implies that $R$ consists of an edge $uv$ such that
$\I_{uv}\setminus \I_S=\{M_2,M_3\}$ and another edge $wx$ such that
$\I_{wx}\setminus \I_S=\{M_1\}$.  Since the third row and column of $P_1$, $P_2$,
and $P_3$ are identical, we see that if any vertex, say $u$, has
weight $v_3$, then the edge in $R$ incident to the vertex must have
either $\I_{uv}\setminus\I_S=\I_R\setminus\I_S$ or
$\I_{uv}\setminus\I_S=\emptyset$.  Neither of these cases occur, so
$u$, $v$, $w$, and $x$ each must have weight $v_1$ or $v_2$.  As in
Subcase~\ref{subcase:setup-for-s-3} in Section~\ref{sec:ltimes}, since $P_1[1,2]$ and $P_2[1,2]=P_3[1,2]$ are
complementary, $\setcard{R}=2$, and $\setcard{\bar \I_S}=1$, we can
apply Corollary~\ref{cor:s-4-to-3} to conclude that $\setcard{S}=3$,
$\setcard{G-H}\leq 3$, and $\setcard{G}\leq 8$.   This contradicts the
assumption that $\setcard{G}\geq 9$, so this case does not occur.
\end{subcase}
\begin{subcase}[$\I_R=\{M_2,M_3,M_4\}$]
Since $\I_R\cup\I_T=\M(H)$, properties~P\ref{prop:IS-U-IT-equiv-classes} and P\ref{prop:IS-U-IT-neq-MH} imply that $\I_T=\{M_1,M_2\}$ or $\I_T=\{M_1,M_2,M_4\}$.
In each of these cases, $\setcard{R}\leq1$ and $\setcard{T}\leq1$ by
Theorem~\ref{thm:minimize-rst}\ref{prop:min-r} and
Corollary~\ref{thm:rt-size}, implying that $\setcard{G-H}\leq3$ and
$\setcard{G}\leq 8$.  This contradicts the assumption that
$\setcard{G}\geq 9$, so these cases do not occur.
\end{subcase}
\end{case}
\begin{case}[$\I_S=\{M_1,M_2,M_3\}$ or $\I_S=\{M_1,M_2,M_4\}$]
In each of these cases, by property~P\ref{prop:extra-class-in-IR},
$\I_R=\{M_1,M_2,M_3,M_4\}$, so $T=\emptyset$ and $\setcard{G-H}\leq
2\setcard{R}$ by Corollary~\ref{thm:t-empty}.  In each of these cases,
$\setcard{R}\leq 1$ by Corollary~\ref{thm:rt-size}, so
$\setcard{G-H}\leq 2$ and $\setcard{G}\leq 7$.  This contradicts the
assumption that $\setcard{G}\geq 9$, so these cases do not occur.

\end{case}
\begin{case}[$\I_S=\{M_3,M_4\}$]
By Observation~\ref{obs:wt-in-intersection-col-space}, if $v\in S$,
then 
\begin{equation*}
  \begin{split}
\wt(v)\in\col(M_1)=\col(M_2)=\{\0,\, v_1=(0,0,0,1,1),\,
v_2=(0,1,1,0,0),\, v_3=(0,1,1,1,1),\,\\ v_4=(1,0,0,0,0),\,
v_5=(1,0,0,1,1),\, v_6=(1,1,1,0,0),\, v_7=(1,1,1,1,1)\}.    
  \end{split}
\end{equation*}
 The
rank-preserving tables for $[v_1,v_2,v_3,v_4,v_5,v_6,v_7]$ with
respect to $M_1$ and $M_2$ are, respectively,
\begin{equation*}
P_1=\begin{bmatrix}
0&1&1&1&1&0&0\\
1&1&0&1&0&0&1\\
1&0&1&0&1&0&1\\
1&1&0&0&1&1&0\\
1&0&1&1&0&1&0\\
0&0&0&1&1&1&1\\
0&1&1&0&0&1&1
\end{bmatrix} \quad \text{ and } \quad 
P_2=\begin{bmatrix}
1&0&1&1&0&1&0\\
0&0&0&1&1&1&1\\
1&0&1&0&1&0&1\\
1&1&0&0&1&1&0\\
0&1&1&1&1&0&0\\
1&1&0&1&0&0&1\\
0&1&1&0&0&1&1
\end{bmatrix}.
\end{equation*}
By property~P\ref{prop:extra-class-in-IR}, $\I_R=\{M_1,M_2,M_3,M_4\}$,
so $T=\emptyset$ and $\setcard{G-H}\leq \setcard{S}\leq 2\setcard{R}$
by Corollary~\ref{thm:t-empty}.  By Corollary~\ref{thm:rt-size},
$\setcard{R}\leq 2$.  If $\setcard{R}=1$, then $\setcard{G-H}\leq 2$
and $\setcard{G}\leq 7$, a contradiction.

Suppose that $\setcard{R}=2$.  Let $R=\{uv,wx\}$.
Theorem~\ref{thm:minimize-rst}\ref{prop:min-r} implies that $R$
consists of an edge $e_1$ such that $\I_{e_1}\setminus \I_S=\{M_1\}$
and another edge $e_2$ such that $\I_{e_2}\setminus\I_S=\{M_2\}$.
Since the third, fourth, and seventh rows and columns of $P_1$ and
$P_2$ are identical, we see that if any vertex, say a vertex in $e_1$,
has weight $v_3$, $v_4$, or $v_7$, then the edge in $R$ incident to
the vertex must have either $\I_{e_1}\setminus\I_S=\I_R\setminus\I_S$
or $\I_{e_1}\setminus \I_S=\emptyset$.  Neither of these cases occur,
so each of the vertices in $S$ must have weight $v_1$, $v_2$, $v_5$,
or $v_6$.  Note also that $P_1[1,2,5,6]$ and $P_2[1,2,5,6]$ are
complementary.  However, we cannot proceed as before and apply
Corollary~\ref{cor:s-4-to-3} since $\setcard{\bar \I_S}=2$.

If there are vertices $a$ and $b$, one incident to each edge of $R$, such
that $\I_{\{a,b\}}=\I_S$, then we can apply
Lemma~\ref{lem:s-4-to-3-general} and conclude that $\setcard{S}=3$,
$\setcard{G-H}\leq 3$, and $\setcard{G}\leq 8$, a contradiction.

Suppose that $\setcard{S}=4$ and there are not two vertices $a$ and
$b$ in $R$ such that $a$ is incident to one edge, $b$ is incident to
the other edge, and $\I_{\{a,b\}}=\I_S=\{M_3,M_4\}$.  By relabeling,
if necessary, we then have $\I_u=\{M_3\}$, $\I_v=\{M_4\}$,
$\I_w=\emptyset$, and $\I_x=\emptyset$.  Recall also that for any
vertex $a\in S$, $\wt(a)\in\{v_1,v_2,v_5,v_6\}$.  Notice that if a
vertex $a$ has weight $\wt(a)=v_1$, then $\I_a=\emptyset$, so
$\wt(u)\neq v_1$ and $\wt(v)\neq v_1$.  Moreover, $\wt(u)\in\col(M_4)$
while $v_2,v_5\not\in\col(M_4)$.  Thus $\wt(u)=v_6$.  Also
$\wt(v)\in\col(M_3)$ and $v_5,v_6\not\in\col(M_3)$, so $\wt(v)=v_2$.
Since $\wt(w),\wt(x)\in\col(M_i)$ for all $i$, $\wt(w)=\wt(x)=v_1$.

Since $\setcard{R}=2$, either $\I_{uv}\setminus\I_S=\{M_1\}$ and
$\I_{wx}\setminus\I_S=\{M_2\}$, or $\I_{uv}\setminus\I_S=\{M_2\}$ and $\I_{wx}\setminus\I_S=\{M_1\}$.

Suppose that $\I_{uv}\setminus\I_S=\{M_1\}$ and
$\I_{wx}\setminus\I_S=\{M_2\}$.  Since $M_2\in\I_{wx}$, $\wt(wx)=0$,
which implies that $M_3\in\I_{wx}$.  Either $M_2\in\I_{vw}$ or
$M_2\not\in\I_{vw}$.
\begin{enumerate}[{\listindent}I.]
\item $M_2\in\I_{vw}$.  Let $R'=\{uv,vw\}$.
\item $M_2\not\in\I_{vw}$.  Then $\wt(vw)=0$, so $M_1\in\I_{vw}$.  Let
  $R'=\{vw,wx\}$.
\end{enumerate}
In either case, $\I_{R'}=\M(H)$, so $G\not\in\F_4(\fh)$ by Corollary~\ref{cor:proper-subset-neq-MH}.  This is a contradiction.

Suppose that $\I_{uv}\setminus\I_S=\{M_2\}$ and
$\I_{wx}\setminus\I_S=\{M_1\}$.  Either $M_1\in\I_{vw}$ or
$M_1\not\in\I_{vw}$.
\begin{enumerate}[{\listindent}I.]
\item $M_1\in\I_{vw}$.  Let $R'=\{uv,vw\}$.
\item $M_1\not\in\I_{vw}$.  Then $\wt(vw)=1$, so $M_2\in\I_{vw}$.
  Also, as can easily be checked, $M_3\in\I_{vw}$.  Let
  $R'=\{vw,wx\}$.
\end{enumerate}
In either case, $\I_{R'}=\M(H)$, so $G\not\in\F_4(\fh)$ by Corollary~\ref{cor:proper-subset-neq-MH}.  This is a contradiction.

Therefore $\setcard{S}\neq 4$, so $\setcard{G-H}\leq \setcard{S}\leq
3$ and $\setcard{G}\leq 8$.  This contradicts the assumption that
$\setcard{G}\geq 9$, so this case does not occur.
\end{case}

For every possible value of $\I_S$, we have reached a contradiction.
Thus $\setcard{G-H}\geq 4$ is impossible, so $\setcard{G-H}\leq 3$ and
$\setcard{G}\leq 8$.
\end{proof}

\section{$H=P_4$}
\label{sec:p4-lemma}

\begin{lemma}\label{lem:p4}
  If $G\in\F_4(P_4)$, then $\setcard{G}\leq 8$.
\end{lemma}
\begin{proof}
Suppose that $G\in \F_4(P_4)$ and $\setcard{G}\geq 8$ (i.e.,
$\setcard{G-H}\geq 4$).  Then $G-H$ has no vertices with zero weight
by Corollary~\ref{thm:zero-vector}.  Assume that $(R,S,T)$ is an
optimal triple for $G$ and $P_4$.  Let
$\M(P_4)=\{M_1,M_2,M_3,M_4,M_5\}$ and
$\C(P_4)=\{C_1=\{M_1,M_2\},C_2=\{M_3,M_4\},C_3=\{M_5\}\}$ be as in
Proposition~\ref{thm:M-calc}\ref{thm:M-P4}.   By
properties~P\ref{prop:IS-U-IT-equiv-classes} and
P\ref{prop:IS-U-IT-neq-MH}, $\I_S\in\{\emptyset, C_1,C_2,C_3,C_1\cup
C_2, C_1\cup C_3, C_2\cup C_3\}$, so we
have the following cases.

\begin{case}[$\I_S=\emptyset$]
By Observation~\ref{obs:wt-in-intersection-col-space}, if $v\in S$,
then $$\wt(v)\in\bigcap_{i=1}^5\col(M_i)=\{\0,\, v_1=(1,0,0,1)^T\}.$$ The
rank-preserving tables for $[v_1]$ with respect to $M_1$, $M_2$, $M_3$, $M_4$, and $M_5$ are, respectively,
\begin{align*}
P_1=
\begin{bmatrix}
  1
\end{bmatrix}, \quad
P_2=
\begin{bmatrix}
  0
\end{bmatrix}, \quad
P_3=
\begin{bmatrix}
  1
\end{bmatrix}, \quad
P_4=
\begin{bmatrix}
  0
\end{bmatrix}, \quad \text{ and }
P_5=
\begin{bmatrix}
  1
\end{bmatrix}.
\end{align*}
Since $P_1=P_3=P_5$ and $P_2=P_4$ are complementary, an edge $uv$ in
$R$ has either $\I_{uv}=\{M_1,M_3,M_5\}$ or $\I_{uv}=\{M_2,M_4\}$.  Thus
$\I_R\in\{\{M_1,M_2,M_3,M_4,M_5\},\, \{M_1,M_3,M_5\},\, \{M_2,M_4\}\}$.  By property~P\ref{prop:extra-class-in-IR}, $\I_R\neq \{M_2,M_4\}$.  Therefore we have
the following cases.

\begin{subcase}[$\I_R=\{M_1,M_2,M_3,M_4,M_5 \}$]
Since $\I_R=\M(P_4)$, Corollary~\ref{thm:t-empty} implies that
$T=\emptyset$ and $\setcard{G-H}\leq 2\setcard{R}$.  By
Theorem~\ref{thm:minimize-rst}\ref{prop:min-r}, $\setcard{R}\leq 2$,
so $\setcard{G-H}\leq 4$ and $\setcard{G}\leq 8$.

\end{subcase}
\begin{subcase}[$\I_R=\{M_1,M_3,M_5\}$]
Since $\I_R\cup\I_T=\M(H)$,
properties~P\ref{prop:IS-U-IT-equiv-classes} and
P\ref{prop:IS-U-IT-neq-MH} imply that $\I_T=\{M_1,M_2,M_3,M_4\}$.  By
Corollary~\ref{thm:rt-size}, $\setcard{R}\leq 1$ and $\setcard{T}\leq
2$.  By Theorem~\ref{thm:gen-bounds}, $\setcard{G-H}\leq 4$, so
$\setcard{G}\leq 8$.
\end{subcase}
\end{case}
\begin{case}[$\I_S=\{M_1,M_2\}$]
By Observation~\ref{obs:wt-in-intersection-col-space}, if $v\in S$,
then 
\begin{equation*}
  \begin{split}
 \wt(v)\in\col(M_3)\cap\col(M_4)\cap\col(M_5)=\{\0,\, v_1=(0,1,0,1)^T,\, v_2=(1,0,0,1)^T,\,\\ v_3=(1,1,0,0)^T\}. 
  \end{split}
\end{equation*}
The rank-preserving
tables for $[v_1,v_2,v_3]$ with respect to $M_3$, $M_4$, and $M_5$ are,
respectively,
\begin{align*}
P_3=
  \begin{bmatrix}
    0&0&0\\
    0&1&1\\
    0&1&1
\end{bmatrix}, \quad
P_4=
  \begin{bmatrix}
1&1&0\\
1&0&1\\
0&1&1\\
\end{bmatrix}, \quad \text{ and }
P_5=
  \begin{bmatrix}
0&0&0\\
0&1&1\\
0&1&1\\
\end{bmatrix}.
\end{align*}
Since $P_3=P_5$, an edge in $R$ is rank-preserving for $M_3$ if and
only if it is also rank-preserving for $M_5$.  Thus
$\I_R\in\{\{M_1,M_2,M_3,M_4,M_5\},\, \{M_1,M_2,M_3,M_5\},\,
\{M_1,M_2,M_4\}\}$.  By property~P\ref{prop:extra-class-in-IR},
$\I_R\neq \{M_1,M_2,M_4\}$.  Therefore we have the following cases.
\begin{subcase}[$\I_R=\{M_1,M_2,M_3,M_4,M_5 \}$]
Since $\I_R=\M(P_4)$, Corollary~\ref{thm:t-empty} implies that
$T=\emptyset$ and $\setcard{G-H}\leq 2\setcard{R}$.  By
Theorem~\ref{thm:minimize-rst}\ref{prop:min-r}, $\setcard{R}\leq 2$,
so $\setcard{G-H}\leq 4$ and $\setcard{G}\leq 8$.
\end{subcase}
\begin{subcase}[$\I_R=\{M_1,M_2,M_3,M_5\}$]
Since $\I_R\cup\I_T=\M(H)$,
properties~P\ref{prop:IS-U-IT-equiv-classes} and
P\ref{prop:IS-U-IT-neq-MH} imply that $\I_T=\{M_1,M_2,M_3,M_4\}$ or
$\I_T=\{M_3,M_4\}$.  In each of these cases, $\setcard{R}\leq 1$ and
$\setcard{T}\leq1$ by Theorem~\ref{thm:minimize-rst}\ref{prop:min-r}
and Corollary~\ref{thm:rt-size}, implying that $\setcard{G-H}\leq3$
and $\setcard{G}\leq 7$.  This contradicts the assumption that
$\setcard{G}\geq 8$, so these cases do not occur.
\end{subcase}

\end{case}
\begin{case}[$\I_S=\{M_3,M_4\}$]
By Observation~\ref{obs:wt-in-intersection-col-space}, if $v\in S$,
then 
\begin{equation*}
  \begin{split}
\wt(v)\in\col(M_1)\cap\col(M_2)\cap\col(M_5)=\{\0,\,
v_1=(0,0,1,1)^T,\, v_2=(1,0,0,1)^T,\,\\ v_3=(1,0,1,0)^T\}.    
  \end{split}
\end{equation*}
The rank-preserving tables for $[v_1,v_2,v_3]$ with respect to $M_1$,
$M_2$, and $M_5$ are, respectively,
\begin{align*}
P_1=\begin{bmatrix}
1&1&0\\
1&1&0\\
0&0&0\\
\end{bmatrix},\quad
P_2=\begin{bmatrix}
1&1&0\\
1&0&1\\
0&1&1\\
\end{bmatrix},\quad \text{ and }
P_5=\begin{bmatrix}
1&1&0\\
1&1&0\\
0&0&0\\
\end{bmatrix}.
\end{align*}
Since $P_1=P_5$, an edge in $R$ is rank-preserving for $M_1$ if and
only if it is also rank-preserving for $M_5$.  Thus
$\I_R\in\{\{M_1,M_2,M_3,M_4,M_5\},\, \{M_1,M_3,M_4,M_5\},\,
\{M_2,M_3,M_4\}\}$.  By property~P\ref{prop:extra-class-in-IR},
$\I_R\neq \{M_2,M_3,M_4\}$.  Therefore we have the following cases.

\begin{subcase}[$\I_R=\{M_1,M_2,M_3,M_4,M_5 \}$]
Since $\I_R=\M(P_4)$, Corollary~\ref{thm:t-empty} implies that
$T=\emptyset$ and $\setcard{G-H}\leq 2\setcard{R}$.  By
Theorem~\ref{thm:minimize-rst}\ref{prop:min-r}, $\setcard{R}\leq 2$,
so $\setcard{G-H}\leq 4$ and $\setcard{G}\leq 8$.
\end{subcase}
\begin{subcase}[$\I_R=\{M_1,M_3,M_4,M_5\}$]
Since $\I_R\cup\I_T=\M(H)$,
properties~P\ref{prop:IS-U-IT-equiv-classes} and
P\ref{prop:IS-U-IT-neq-MH} imply that $\I_T=\{M_1,M_2,M_3,M_4\}$ or
$\I_T=\{M_1,M_2\}$.  In each of these cases, $\setcard{R}\leq 1$ and
$\setcard{T}\leq1$ by Theorem~\ref{thm:minimize-rst}\ref{prop:min-r}
and Corollary~\ref{thm:rt-size}, implying that $\setcard{G-H}\leq3$
and $\setcard{G}\leq 7$.  This contradicts the assumption that
$\setcard{G}\geq 8$, so these cases do not occur.
\end{subcase}

\end{case}
\begin{case}[$\I_S=\{M_5\}$]
By Observation~\ref{obs:wt-in-intersection-col-space}, if $v\in S$,
then 
\begin{equation*}
  \begin{split}
\wt(v)\in\col(M_1)\cap\col(M_2)\cap\col(M_3)\cap\col(M_4)=\{\0,\,
v_1=(0,1,1,1)^T,\, v_2=(1,0,0,1)^T,\,\\ v_3=(1,1,1,0)^T\}.    
  \end{split}
\end{equation*}
The rank-preserving tables for $[v_1,v_2,v_3]$ with respect to $M_1$,
$M_2$, $M_3$, and $M_4$ are, respectively,
\begin{align*}
  P_1=
  \begin{bmatrix}
    1&0&1\\
    0&1&1\\
    1&1&0
\end{bmatrix},\quad
P_2=
\begin{bmatrix}
  1&0&1\\
  0&0&0\\
  1&0&1
\end{bmatrix},\quad
P_3=
\begin{bmatrix}
  0&1&1\\
  1&1&0\\
  1&0&1
\end{bmatrix},\quad \text{ and }
P_4=
\begin{bmatrix}
  1&0&1\\
  0&0&0\\
  1&0&1
\end{bmatrix}.
\end{align*}
Since $P_2=P_4$, an edge in $R$ is rank-preserving for $M_2$ if and
only if it is also rank-preserving for $M_4$.  Thus
\begin{equation*}
\begin{split}\I_R\in\{\{M_1,M_2,M_3,M_4,M_5\},\, 
  \{M_1,M_2,M_4,M_5\},\, \{M_2,M_3,M_4,M_5\},\, \{M_2,M_4,M_5\}, \\
  \{M_1,M_3,M_5\},\, \{M_1,M_5\},\, \{M_3,M_5\}\}.
\end{split}
\end{equation*}
By property~P\ref{prop:extra-class-in-IR}, $\I_R\not\in \{
\{M_2,M_4,M_5\},\, \{M_1,M_3,M_5\},\, \{M_1,M_5\},\, \{M_3,M_5\}\}$.
Therefore we have the following cases.

\begin{subcase}[$\I_R=\{M_1,M_2,M_3,M_4,M_5\}$]
Since $\I_R=\M(P_4)$, Corollary~\ref{thm:t-empty} implies that
$T=\emptyset$ and $\setcard{G-H}\leq \setcard{S}\leq 2\setcard{R}$.  By
Theorem~\ref{thm:minimize-rst}\ref{prop:min-r}, $\setcard{R}\leq
3$.  If $\setcard{R}\leq 2$, then $\setcard{G-H}\leq 4$ and
$\setcard{G}\leq 8$.

Suppose that $\setcard{R}=3$.  Then
Theorem~\ref{thm:minimize-rst}\ref{prop:min-r} implies that $R$
consists of three edges $e_1$, $e_2$, and $e_3$ such that
$\I_{e_1}\setminus\I_S=\{M_1\}$, $\I_{e_2}\setminus\I_S=\{M_2,M_4\}$,
and $\I_{e_3}\setminus\I_S=\{M_3\}$.

Since the first row and column of $P_1$ and $P_2$ are the same, if an
edge $e\in R$ is incident to a vertex of weight $v_1$, then either
$\{M_1,M_2\}\subseteq \I_e$ or $\{M_1,M_2\}\subseteq\M(H)\setminus
\I_e$.  Therefore $e_1$ and $e_2$ are not incident to vertices with weight
$v_1$.  Since the third row and column of $P_2$ and $P_3$ are the
same, if an edge $e\in R$ is incident to a vertex of weight $v_3$, then either
$\{M_2,M_3\}\subseteq \I_e$ or $\{M_2,M_3\}\subseteq\M(H)\setminus
\I_e$.  Therefore $e_2$ and $e_3$ are not incident to vertices with weight
$v_3$.  Since $P_1[2]=P_3[2]$, if both vertices incident to an edge $e\in R$ have
weight $v_2$, then $\{M_1,M_3\}\subseteq \I_e$ or
$\{M_1,M_3\}\subseteq\M(H)\setminus \I_e$.  Therefore $e_1$ and $e_3$
each are incident to at least one vertex that does not have weight $v_2$.

Therefore $e_1$ must be incident to vertices with weights $v_2$ and
$v_3$ (implying that \mbox{$\wt(e_1)=0$} since $\I_{e_1}=\{M_1\}$) or
incident to vertices with weights $v_3$ and $v_3$ (implying that
$\wt(e_1)=1$).  Each vertex incident to $e_2$ must have weight $v_2$,
which implies that $\wt(e_2)=1$.  The edge $e_3$ must be incident to
vertices with weights $v_1$ and $v_1$ (implying that $\wt(e_3)=1$) or
incident to vertices with weights $v_1$ and $v_2$ (implying that $\wt(e_3)=0$).

Therefore there are at least three vertices $u$, $v$, and $w$ in $S$
such that $u$ is incident to $e_3$, $v$ is incident to $e_2$, $w$ is
incident to $e_1$, $\wt(u)=v_1$,
$\wt(v)=v_2$, and $\wt(w)=v_3$.  Let $R'=\{uv,vw\}$.  Note that since
$v_1,v_3\not\in\col(M_5)$, $M_5\in\I_{R'}$. Suppose that
$\setcard{S}\geq5$.   Then the vertices in $R'\cup e$ for any edge $e\in R$
form a proper subset of $S$.  We now have the following possibilities
for $\I_{R'}$.
\begin{enumerate}[{\listindent}I.]
\item $M_1\not\in\I_{R'}$.  Then $\wt(vw)=1$, which implies that
  $M_2,M_3,M_4\in\I_{R'}$.  Since $\I_{R'\cup e_1}=\M(H)$,
  $G\not\in\F_4(P_4)$ by Corollary~\ref{cor:proper-subset-neq-MH},
  which is a contradiction.
\item $M_2,M_4\not\in\I_{R'}$.  Then $\wt(uv)=0$, which implies that
  $M_3\in\I_{R'}$.  Also \mbox{$\wt(vw)=0$}, which implies that $M_1\in\I_{R'}$.
  Since $\I_{R'\cup e_2}=\M(H)$, $G\not\in\F_4(P_4)$ by
  Corollary~\ref{cor:proper-subset-neq-MH}, which is a contradiction.
\item $M_3\not\in\I_{R'}$.  Then $\wt(uv)=1$, which implies that
  $M_1,M_2,M_4\in\I_{R'}$.  Since $\I_{R'\cup e_3}=\M(H)$,
  $G\not\in\F_4(P_4)$ by Corollary~\ref{cor:proper-subset-neq-MH},
  which is a contradiction.
\item $\I_{R'}=\M(P_4)$.  Since the vertices in $R'$ are a proper
  subset of the vertices in $S$, $G\not\in\F_4(P_4)$ by
Corollary~\ref{cor:proper-subset-neq-MH}, which is a contradiction.
\end{enumerate}
Since each case leads to a contradiction, our assumption that
$\setcard{S}\geq 5$ must be false.  Therefore $\setcard{S}\leq4$, so
$\setcard{G-H}\leq4$ and $\setcard{G}\leq8$.

\end{subcase}
\begin{subcase}[$\I_R=\{M_1,M_2,M_4,M_5\}$]
Since $\I_R\cup\I_T=\M(H)$,
properties~P\ref{prop:IS-U-IT-equiv-classes} and
P\ref{prop:IS-U-IT-neq-MH} imply that $\I_T=\{M_3,M_4,M_5\}$ or
$\I_T=\{M_3,M_4\}$.  In each of these cases, $\setcard{T}\leq 1$ by
Corollary~\ref{thm:rt-size}.  Since $M_2\in\I_{uv}$ if and only if
$M_4\in\I_{uv}$ for each edge $uv\in R$, $\setcard{R}\leq 2$ by
Theorem~\ref{thm:minimize-rst}\ref{prop:min-r}.  If $\setcard{R}=1$,
then $\setcard{G-H}\leq 3$, a contradiction.  If $\setcard{R}=2$, then
Theorem~\ref{thm:minimize-rst}\ref{prop:min-r} implies that $R$
consists of an edge $uv$ such that $\I_{uv}\setminus
\I_S=\{M_2,M_4\}$  and another
edge $wx$ such that $\I_{wx}\setminus \I_S=\{M_1\}$.  Since the first row and column
of $P_1$, $P_2$, and $P_4$ are identical, we see that if any vertex,
say $u$, has weight $v_1$, then the edge in $R$ incident to the vertex
must have either $\I_{uv}\setminus\I_S=\I_R\setminus\I_S$ or
$\I_{uv}\setminus\I_S=\emptyset$.  Neither of these cases occur, so
$u$, $v$, $w$, and $x$ each must have weight $v_2$ or $v_3$.  As in
Subcase~\ref{subcase:setup-for-s-3} in Section~\ref{sec:ltimes}, since $P_1[2,3]$ and $P_2[2,3]=P_4[2,3]$ are
complementary, $\setcard{R}=2$, and $\setcard{\bar \I_S}=1$, we can
apply Corollary~\ref{cor:s-4-to-3} to conclude that $\setcard{S}=3$,
implying that $\setcard{G-H}\leq 4$ and $\setcard{G}\leq 8$.

\end{subcase}
\begin{subcase}[$\I_R=\{M_2,M_3,M_4,M_5\}$]
Since $\I_R\cup\I_T=\M(H)$,
properties~P\ref{prop:IS-U-IT-equiv-classes} and
P\ref{prop:IS-U-IT-neq-MH} imply that $\I_T=\{M_1,M_2\}$ or
$\I_T=\{M_1,M_2,M_5\}$.  In each of these cases, $\setcard{T}\leq 1$
by Corollary~\ref{thm:rt-size}.  Since $M_2\in\I_{uv}$ if and only if
$M_4\in\I_{uv}$ for each edge $uv\in R$, $\setcard{R}\leq 2$ by
Theorem~\ref{thm:minimize-rst}\ref{prop:min-r}.  If $\setcard{R}=1$,
then $\setcard{G-H}\leq 3$, a contradiction.  If $\setcard{R}=2$, then
Theorem~\ref{thm:minimize-rst}\ref{prop:min-r} implies that $R$
consists of an edge $uv$ such that $\I_{uv}\setminus \I_S=\{M_2,M_4\}$ and another
edge $wx$ such that $\I_{wx}\setminus \I_S=\{M_3\}$.  Note that the third row and
column of $P_2$, $P_3$, and $P_4$ are identical; as in the previous
case, none of $u$, $v$, $x$, $w$ can have weight $v_3$, so each must
have weight $v_1$ or $v_2$.  Since $P_3[1,2]$ and $P_2[1,2]=P_4[1,2]$
are complementary, $\setcard{R}=2$, and $\setcard{\bar \I_S}=1$, we
can apply Corollary~\ref{cor:s-4-to-3} to conclude that
$\setcard{S}=3$, implying that $\setcard{G-H}\leq 4$ and
$\setcard{G}\leq 8$.

\end{subcase}
\end{case}
\begin{case}[$\I_S=\{M_1,M_2,M_3,M_4\}$, $\I_S=\{M_1,M_2,M_5\}$, or %
  $\I_S=\{M_3,M_4,M_5\}$]
In each of these cases, by property~P\ref{prop:extra-class-in-IR},
$\I_R=\{M_1,M_2,M_3,M_4,M_5\}$, so $T=\emptyset$ and
$\setcard{G-H}\leq 2\setcard{R}$ by Corollary~\ref{thm:t-empty}.
In each of these cases, $\setcard{R}\leq 2$ by
Corollary~\ref{thm:rt-size}, so $\setcard{G-H}\leq 4$ and
$\setcard{G}\leq 8$. \qedhere
\end{case}
\end{proof}

\section{All graphs in $\F_4(\Fld_2)$}\label{sec:graphs-in-F_4}
Combining Lemmas~\ref{lem:3k2-p3Vp3} through \ref{lem:p4} with
Corollary~\ref{cor:F4-in-F3-H}, we have:

\begin{theorem}\label{thm:F_4-vertices}
  All graphs in $\F_4(\Fld_2)$ have 8 or fewer vertices.
\end{theorem}

Theorem 3.1 in \cite{ding-kotlov-minrank-finite} implies that all graphs in
$\F_4(\Fld_2)$ have 25 or fewer vertices.  Because we have
made a much more detailed analysis for the field $\Fld_2$, we have
been able to greatly improve their bound in this single case.  Since
all graphs in $\F_4(\Fld_2)$ have 8 or fewer vertices, we can do an
exhaustive search for all the graphs.  In
Appendix~\ref{sec:magma-programs}, we list a few Magma functions
sufficient to implement this search.  These functions use the graph
generation program ``geng'' distributed with Brendan McKay's Nauty
program \cite[Version 2.2]{mckay-nauty}.  This exhaustive search
results in the 62 graphs displayed at the end of this section.  Thus,
recalling Observation~\ref{obs:Gk-iff-Fkplus1-forbidden}, we have:

\begin{theorem}\label{thm:F_4}
  $\F_4(\Fld_2)$ consists of the 62 graphs listed at the end of this
  section.  For any graph $G$, $\mr(\Fld_2,G)\leq3$ if and only if no
  graph in $\F_4(\Fld_2)$ is induced in $G$.
\end{theorem}

In the listing of the graphs in $\F_4(\Fld_2)$ that follows, the
graphs are sorted by number of vertices.  We have also tried to group
similarly structured graphs together.  Each graph is identified with a
number and a graph6 code.  The graph6 code is a compact representation
of the adjacency matrix (and thus the zero/nonzero pattern of the
matrices associated with the graph).  The specification of the graph6
code is distributed with Nauty and can also be found on the Nauty website.

We now proceed with the listing of all 62 graphs in $\F_4(\Fld_2)$.

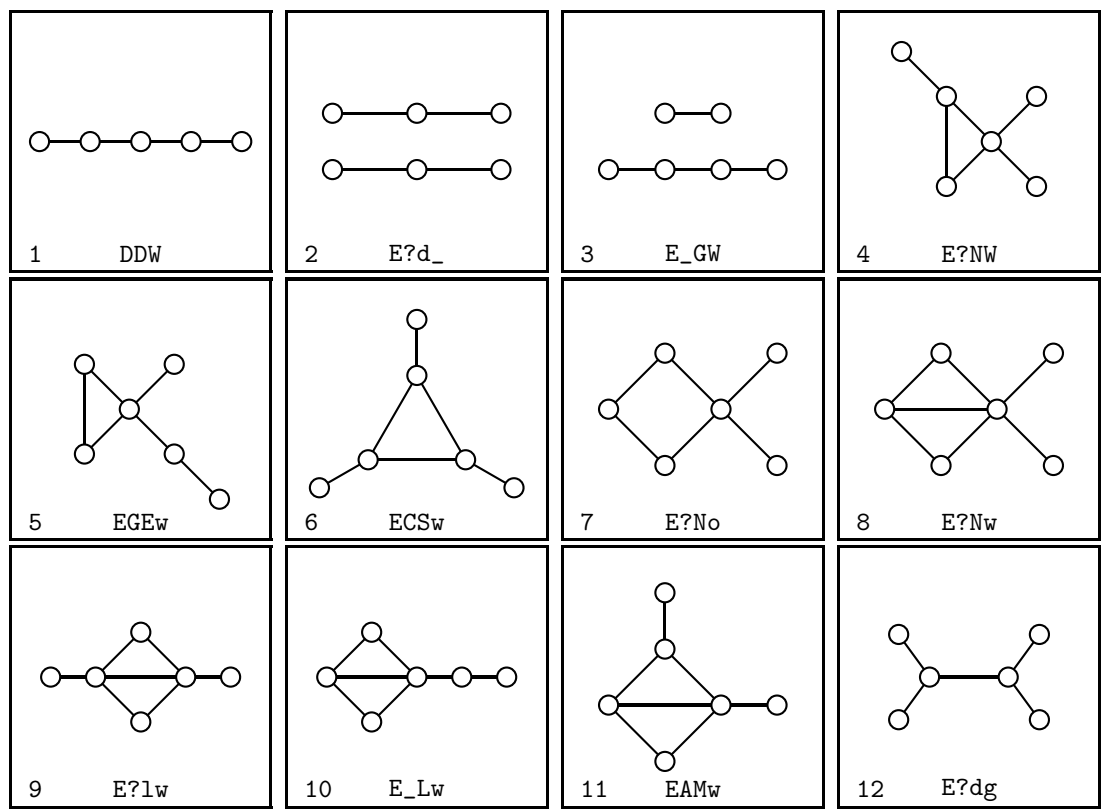
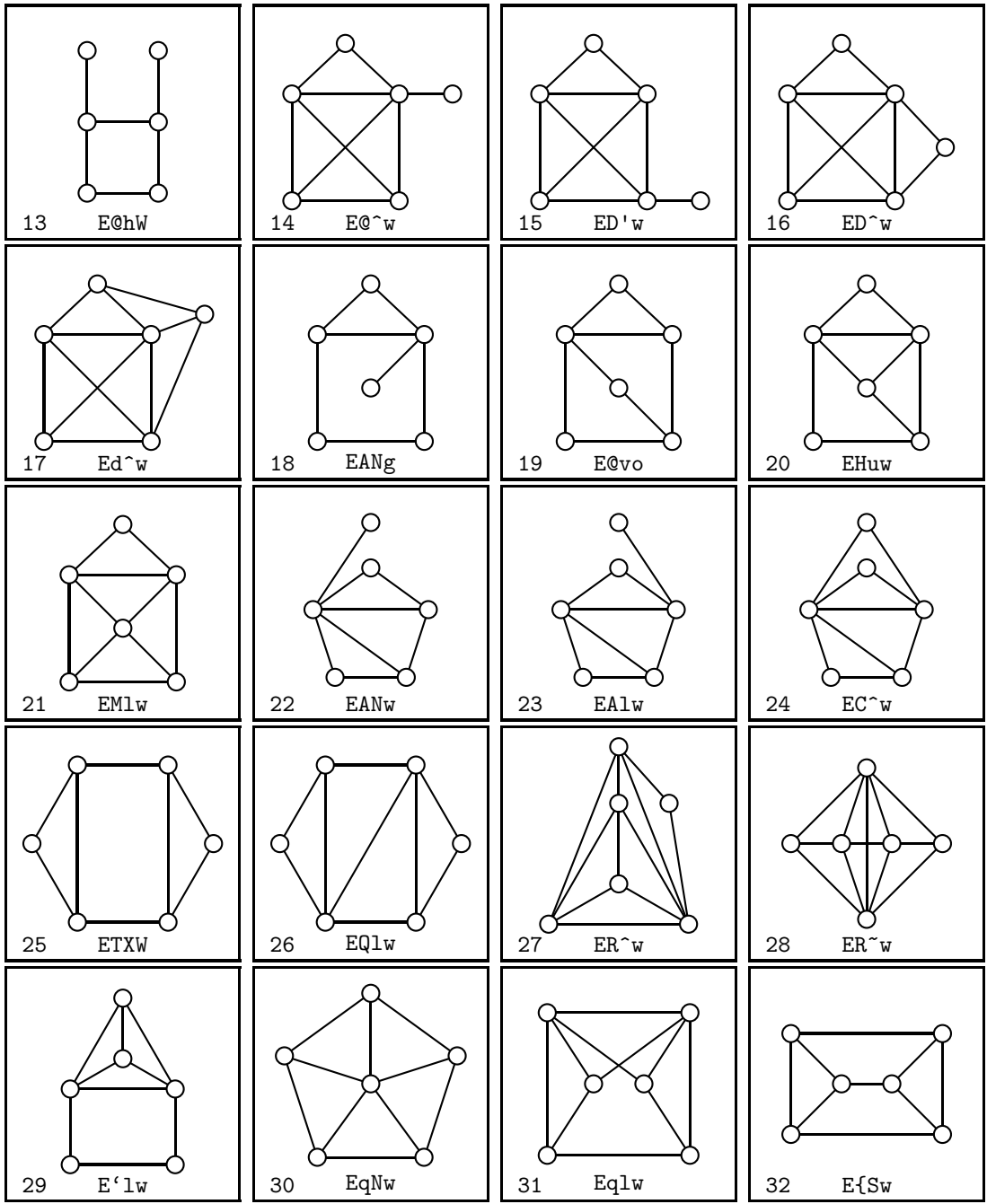
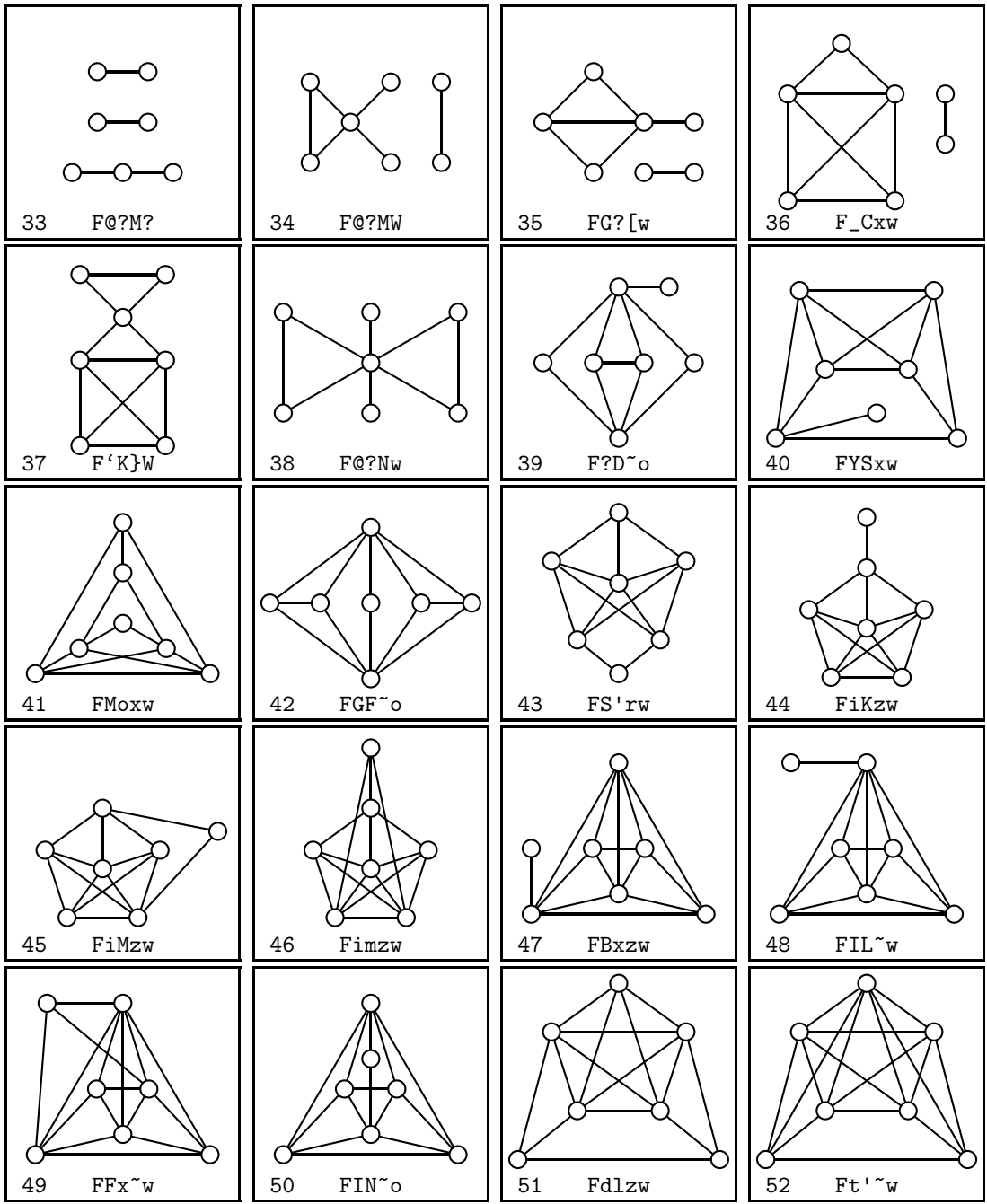
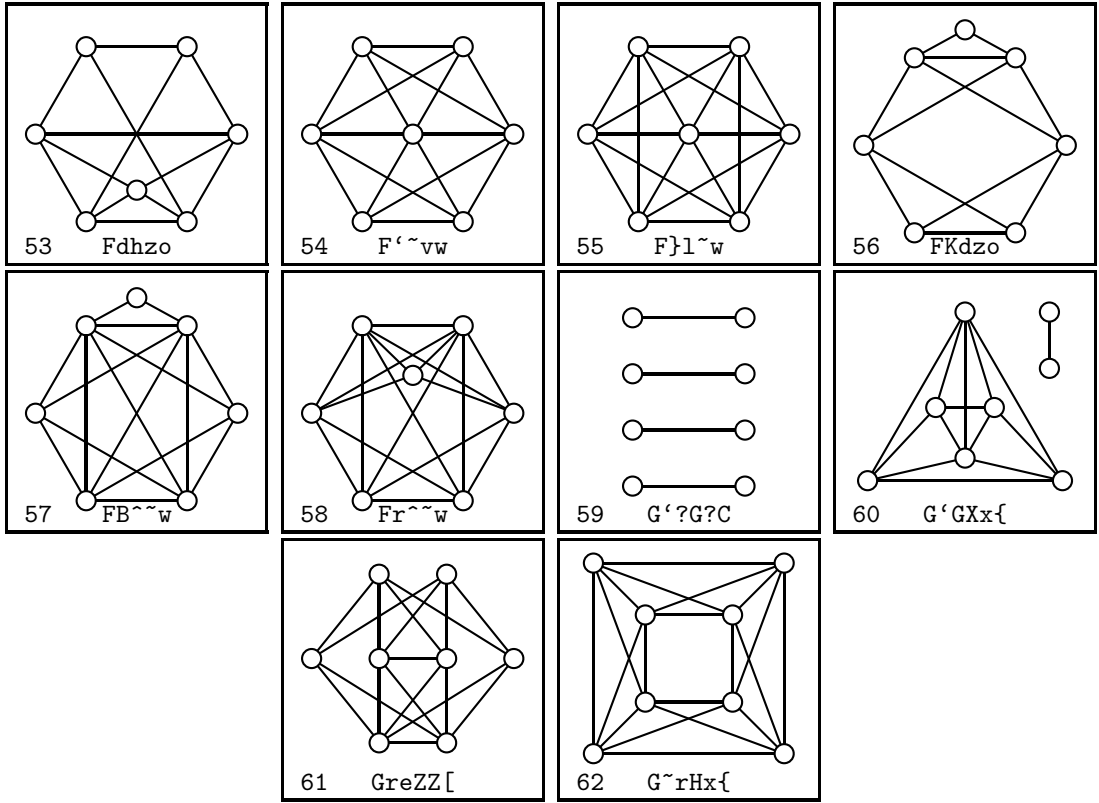
\begin{tabbing} 
\input graphs/graph-1.tex \= \input graphs/graph-2.tex \= \input graphs/graph-3.tex \= \input graphs/graph-4.tex \\
\input graphs/graph-5.tex \> \input graphs/graph-6.tex \> \input graphs/graph-7.tex \> \input graphs/graph-8.tex \\
\input graphs/graph-9.tex \> \input graphs/graph-10.tex \> \input graphs/graph-11.tex \> \input graphs/graph-12.tex \\
\input graphs/graph-13.tex \> \input graphs/graph-14.tex \> \input graphs/graph-15.tex \> \input graphs/graph-16.tex \\
\input graphs/graph-17.tex \> \input graphs/graph-18.tex \> \input graphs/graph-19.tex \> \input graphs/graph-20.tex \\
\input graphs/graph-21.tex \> \input graphs/graph-22.tex \> \input graphs/graph-23.tex \> \input graphs/graph-24.tex \\
\input graphs/graph-25.tex \> \input graphs/graph-26.tex \> \input graphs/graph-27.tex \> \input graphs/graph-28.tex \\
\input graphs/graph-29.tex \> \input graphs/graph-30.tex \> \input graphs/graph-31.tex \> \input graphs/graph-32.tex \\
\input graphs/graph-33.tex \> \input graphs/graph-34.tex \> \input graphs/graph-35.tex \> \input graphs/graph-36.tex \\
\input graphs/graph-37.tex \> \input graphs/graph-38.tex \> \input graphs/graph-39.tex \> \input graphs/graph-40.tex \\
\input graphs/graph-41.tex \> \input graphs/graph-42.tex \> \input graphs/graph-43.tex \> \input graphs/graph-44.tex \\
\input graphs/graph-45.tex \> \input graphs/graph-46.tex \> \input graphs/graph-47.tex \> \input graphs/graph-48.tex \\
\input graphs/graph-49.tex \> \input graphs/graph-50.tex \> \input graphs/graph-51.tex \> \input graphs/graph-52.tex \\
\input graphs/graph-53.tex \> \input graphs/graph-54.tex \> \input graphs/graph-55.tex \> \input graphs/graph-56.tex \\
\input graphs/graph-57.tex \> \input graphs/graph-58.tex \> \input graphs/graph-59.tex \> \input graphs/graph-60.tex \\
\> \input graphs/graph-61.tex \> \input graphs/graph-62.tex 
\end{tabbing}

\section{Graphs in $\F_4(F)$ for other fields}\label{sec:other-fields}
Many of the graphs in $\F_4(\Fld_2)$ are also in $\F_4(F)$ for any
field $F$.  This is the case with most of the disconnected graphs and
the connected graphs with a cut vertex in the table.

We need the following elementary facts \cite{barrett-vdHL-minrank2-infinite}.  

\begin{observation}\label{obs:other-fields}
  For any field $F$
  \begin{enumerate}
  \item $\mr(F,K_n)=1$ for $n\geq 2$;
    $\mr(F,K_{2,3})=\mr(F,\input{graphs/graph39_outside_pendant.tex})=2$;
    $\mr(F,\ltimes)=\mr(F,\dart)=3$.
      \item\label{obs-part:specific-graphs} $K_2\in\F_1(F)$; $\ltimes,\dart\in\F_3(F)$.
      \item\label{obs-part:disconnected} If $G=\cup_{i=1}^k G_i$, then $\mr(F,G)=\sum_{i=1}^k \mr(F,G_i)$.
  \end{enumerate}
\end{observation}
We will also need
\begin{theorem}[{\cite{fiedler-tridiag-matrices,bento-duarte-tridiag-matrices}}]\label{thm:pn-rank-n-1}
  Let $F$ be any field and let $G$ be a graph on $n$ vertices.  Then
  $\mr(F,G)=n-1$ if and only if $G=P_n$.
\end{theorem}
A stronger result was proved by Fiedler over $\mathbb{R}$ \cite{fiedler-tridiag-matrices}
and his result was extended to any field, with some exceptions for $\Fld_3$, by Bento and Duarte \cite{bento-duarte-tridiag-matrices}.

\begin{corollary}\label{cor:path-mr}
  For any field $F$, $\mr(F,P_n)=n-1$ and $P_n\in\F_{n-1}(F)$.
\end{corollary}

We will also utilize the following

\begin{proposition}\label{prop:mr3-class}
  Let $\mathcal{E}=\{\fh, G_1=\input{graphs/graph40_pendant.tex},
  G_2=\input{graphs/graph44_pendant.tex}, P_3\lor P_3\}$ ($G_1$ is
  graph 40 minus the pendant vertex and $G_2$ is graph 44 minus the
  pendant vertex).  Then for each $G\in\mathcal{E}$, $\mr(\Fld_2,G)=3$ and $\mr(F,G)=2$ for any $F\neq \Fld_2$.  Moreover, $\fh,P_3\lor P_3\in\F_3(\Fld_2)$.
\end{proposition}

\begin{proof}
  We already verified the first claim for the full house in the introduction.
  Taking complements of the others we find that
  $\comp{G_1}=2P_3$, $\comp{G_2}=P_3\cup K_2\cup K_1$, and
  $\comp{(P_3\lor P_3)}=2K_2\cup2K_1$.  By Theorems~6 and 7 in
  \cite{barrett-vdHL-minrank2-infinite} and Theorems~11 and 15 in
  \cite{barrett-vdHL-minrank2-finite},
  $\mr(F,G_1)=\mr(F,G_2)=\mr(F,P_3\lor P_3)=2$ for $F\neq \Fld_2$, while
  $\mr(\Fld_2,G_1)=\mr(\Fld_2,G_2)=\mr(\Fld_2,P_3\lor P_3)=3$.  The final claim
  follows from Theorem~\ref{thm:F3-F2}.
\end{proof}

\subsection{Disconnected graphs}
\begin{proposition}\label{prop:disconnected-sum-forbidden-mr}
  If $F$ is any field and $S_i\in\F_{\mr(S_i)}(F)$, $i=1,\ldots,m$, then
  \begin{equation*}
    \bigcup_{i=1}^m S_i \in \F_{\mr(S_1)+\cdots+\mr(S_m)}(F).
  \end{equation*}
\end{proposition}
\begin{proof}
  This follows immediately from
  Observation~\ref{obs:other-fields}\ref{obs-part:disconnected} and
  the definition of $\F_{k+1}(F)$.
\end{proof}

Applying Observation
\ref{obs:other-fields}\ref{obs-part:specific-graphs},
Corollary~\ref{cor:path-mr}, and
Proposition~\ref{prop:disconnected-sum-forbidden-mr} to the
disconnected graphs 2, 3, 33, 34, 35, and 59 in
Section~\ref{sec:graphs-in-F_4}, we have

\begin{theorem}
  For any field $F$,
  \begin{equation*}
    \{2P_3,P_4\cup K_2,P_3\cup2K_2,\ltimes \cup K_2, \dart \cup K_2, 4K_2\}\subseteq\F_4(F).
  \end{equation*}
\end{theorem}

Graphs 36 and 60 are $\fh\, \cup K_2$ and $(P_3\lor P_3)\cup K_2$.  Since $\fh,
P_3\lor P_3\in\F_3(F)$ if and only if $F=\Fld_2$, graphs 36 and 60 are not in
$\F_4(F)$ for any $F\neq \Fld_2$.

\subsection{Connected graphs with a cut vertex}

First, note that graph~1 in Section~\ref{sec:graphs-in-F_4}, $P_5$, is
in $\F_4(F)$ for any field $F$ by Corollary~\ref{cor:path-mr}.

We now recall a definition and a known result.

\begin{definition}\label{def:vsum}
  Let $G$ and $H$ be graphs on at least two vertices, each having a
  vertex labeled $v$.  Then $G\vsum H$ is the graph obtained from $G\cup
  H$ by identifying the two vertices labeled $v$.  Similarly, if
  $G_1,\ldots,G_k$, $k\geq2$, are graphs on at least two vertices, each with a
  vertex labeled $v$, let $G=G_1\vsum G_2 \vsum \cdots \vsum G_k$ be the
  graph obtained by identifying the vertices labeled $v$ in each of
  the graphs.  We call $G$ the \emph{vertex sum} of the graphs
  $G_1,\ldots,G_k$ at $v$.  Note that $v$ is necessarily a cut vertex of a
  graph constructed in this way and that any graph with a cut vertex
  can be expressed as such a sum with $k\geq2$.
\end{definition}

The following theorem was proved over the real field in
\cite{hsieh-minrank} and \cite{BFH1-minrankpath}.  In
Appendix~\ref{appendix:cut-vertex}, we give an easy proof of part
\ref{thm:cut-vertex-two-graphs} that holds for any field; part
\ref{thm:cut-vertex-many-graphs} then follows by induction.  This same
proof is also a key part of the proof of a more general result on the
inertia set of a graph with a cut vertex (see Theorem~4.2 in
\cite{barrett-loewy-hall-graph-inertia}).

\begin{theorem}[{\cite{hsieh-minrank,BFH1-minrankpath}}]\label{thm:cut-vertex}
  Let $F$ be any field.
  \begin{enumerate}
  \item\label{thm:cut-vertex-two-graphs} If $G$ and $H$ are graphs on
    at least two vertices, each having a vertex labeled $v$, then
    \begin{equation*}
      \mr(F,G \vsum H)=\min\{\mr(F,G)+\mr(F,H),\, \mr(F,G-v)+\mr(F,H-v)+2\}.
    \end{equation*}

  \item\label{thm:cut-vertex-many-graphs} Let $G_1,\ldots,G_k$, $k\geq2$, be
    graphs on at least two vertices, each with a vertex labeled $v$.
    Then
    \begin{equation*}
      \mr(F, G_1\vsum G_2 \vsum \cdots \vsum G_k)=\min\{\sum_{i=1}^k\mr(F,G_i), \sum_{i=1}^k\mr(F,G_i-v)+2\}.
    \end{equation*}

  \end{enumerate}

\end{theorem}

This result reduces the calculation of the minimum rank of any graph
with a cut vertex to a calculation for smaller graphs.

\begin{corollary}\label{cor:vsum-no-more-than-sum}
  $\mr(F,G\vsum H)\leq \mr(F,G)+\mr(F,H)$.
\end{corollary}

We can now establish one criterion for membership in $\F_4(F)$ for any
field.

\begin{theorem}\label{thm:F_4-sufficient-any-field}
  Let $F$ be any field and let $G$ be a graph satisfying all of the
  following
  \begin{enumerate}
  \item $\setcard{G}=6$,
  \item $\mr(F,G)=4$, and
  \item $P_5$ is not induced in $G$.
  \end{enumerate}
  Then $G\in\F_4(F)$.
\end{theorem}

\begin{proof}
  For each vertex $v$ of $G$, $G-v$ is a graph on 5 vertices distinct
  from $P_5$.  By Theorem~\ref{thm:pn-rank-n-1}, $\mr(F,G)<5-1=4$.  By
  Definition~\ref{def:F_k+1}, $G\in\F_4(F)$.
\end{proof}

\begin{proposition}
  Graphs 4, 5, 6, 7, 8, 9, 10, 11, 12, 13, 18, 22, and 23 are all in
  $\F_4(F)$ for any field $F$.
\end{proposition}

\begin{proof}
  Each of these graphs has 6 vertices and $P_5$ is induced in none of
  them. Moreover, each graph is of the form $G \vsum K_2$, where $G\neq
  \fh$ is a graph on 5 vertices.  Let $G \vsum K_2$ be any of these
  graphs.  By Theorem~\ref{thm:cut-vertex},
  Proposition~\ref{prop:no-fh-indep-field}, and Theorem~\ref{thm:F_4},
  \begin{align*}
    \mr(F,G \vsum K_2)&= \min\{\mr(F,G)+\mr(F,K_2),\mr(F,G-v)+\mr(F,K_2-v)+2\}\\
    &=\min\{\mr(\Fld_2,G)+\mr(\Fld_2,K_2),\, \mr(\Fld_2,G-v)+\mr(\Fld_2,K_2-v)+2\}\\
    &=\mr(\Fld_2,G\vsum K_2)=4.
  \end{align*}
  By Theorem~\ref{thm:F_4-sufficient-any-field}, $G\vsum K_2\in\F_4(F)$.
\end{proof}

We note that graphs 14 and 15, which contain the full house, have
minimum rank 3 over any field $F\neq \Fld_2$, so are not in $\F_4(F)$ for
$F\neq \Fld_2$.\bigskip

We now consider in turn graphs 38 and 39.

Graph 38 (\input{graphs/graph38inline.tex}): Applying
Theorem~\ref{thm:cut-vertex}\ref{thm:cut-vertex-many-graphs} with
$k=4$, we have 
\begin{align*}
  \mr(F,\input{graphs/graph38inline.tex}) &= \min\{2\mr(F,K_3)+2\mr(F,K_2),\ 2\mr(F,K_2)+2\mr(F,K_1)+2\}\\
  &= \min \{2+2,2+0+2\} = 4
\end{align*}
Theorem~\ref{thm:cut-vertex} also implies that
$\mr(F,\input{graphs/graph38deleteinline.tex})=3$.  By definition,
$\input{graphs/graph38inline.tex} \in \F_4(F)$.

  \begin{figure}[ht]
    \centering
    \input{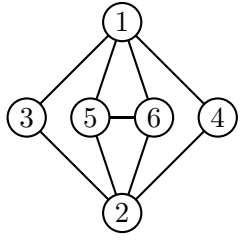}
    \caption{$H=$ graph 39 minus the pendant vertex.}
    \label{fig:g39-v}
  \end{figure}
Graph 39 (\input{graphs/graph39inline.tex}): Let $F$ be any
  field.  Let $H$ be the graph obtained by deleting the pendant vertex
  in graph 39, labeled as in Figure~\ref{fig:g39-v}.
  Since $\ltimes$ is induced in $H$, $\mr(F,H)\geq\mr(F,\ltimes)=3$ by
  Observation~\ref{obs:other-fields}.  Moreover,
  \begin{equation*}
    A=
    \begin{bmatrix}
      0&0&1&1&1&1\\
      0&0&1&1&1&1\\
      1&1&0&0&0&0\\
      1&1&0&0&0&0\\
      1&1&0&0&1&1\\
      1&1&0&0&1&1
    \end{bmatrix}\in\mathcal{S}(F,H)
  \end{equation*}
  and $\rank A=3$.  Therefore $\mr(F,H)=3$.  By
  Theorem~\ref{thm:cut-vertex},
  \begin{align*}
    \mr(F,\textrm{graph 39}) &=
    \min\{\mr(F,H)+\mr(F,K_2),\,\mr(F,\ltimes)+\mr(F,K_1)+2\}\\
    &= \min\{ 3+1,\,3+0+2\}=4.
  \end{align*}
  Any graph obtained by deleting a vertex from graph 39 is one of $H$,
  $\ltimes\cup K_1$, \input{graphs/graph38deleteinline.tex},
  \mbox{\input{graphs/graph39_inside},} or
  \input{graphs/graph39_outside}.  By
  Observation~\ref{obs:other-fields}, $\mr(F,\ltimes\cup K_1)=3$.  We
  just saw that \input{graphs/graph38deleteinline.tex}\, has minimum
  rank 3.  Since $K_{2,3}$ and
  \input{graphs/graph39_outside_pendant.tex} have minimum rank 2 over
  any field by Observation~\ref{obs:other-fields}, the graphs
  \input{graphs/graph39_inside} and \input{graphs/graph39_outside}
  each have minimum rank at most 3 by
  Corollary~\ref{cor:vsum-no-more-than-sum}.  By definition,
  $\textrm{graph 39}\in\F_4(F)$ for every field $F$.

Summarizing,
\begin{proposition}
  Graphs 38 and 39 are in $\F_4(F)$ for any field $F$.
\end{proposition}

The four remaining connected graphs with cut vertices, graphs 40, 44,
47, and 48, in the table do not belong to $\F_4(F)$ for $F\neq \Fld_2$.  Let
$G$ be any of these graphs.  Deleting the pendant vertex in $G$ yields
one of the last three graphs in Proposition~\ref{prop:mr3-class}, so
by that result and Corollary~\ref{cor:vsum-no-more-than-sum}, 
$\mr(F,G)\leq2+1=3$ for $F\neq \Fld_2$.

\subsection{Summary}

We have seen that 6 of the 8 disconnected graphs in
Section~\ref{sec:graphs-in-F_4} are in $\F_4(F)$ for all fields $F$,
while 16 of 22 of the connected graphs with a cut vertex are in
$\F_4(F)$ for all $F$.

We stated in the introduction that even if one is only interested in
the minimum rank problem over $\mathbb{R}$, results obtained over
$\Fld_2$ yield important insights.  We have just observed that of the
30 graphs with vertex connectivity at most one in the list of 62
graphs in $\F_4(\Fld_2)$, 22 of these are also in $\F_4(F)$ for
\emph{any} field.  While the discrepancy is significant, it is also
the case that the amount of overlap is surprising.  The analysis of
the 2-connected graphs in Section~\ref{sec:graphs-in-F_4} seems to be much more
complicated with our present tools.

We have not found all graphs with vertex connectivity less than 2 in
$\F_4(F)$, $F\neq \Fld_2$, by the above methods.  For example, let $F$ be any field with
$\character F\neq 2$.  Then $\mr(F,K_{3,3,3})=3$ and $\mr(F,K_{3,3,2})=2$
\cite{barrett-vdHL-minrank2-infinite}.  Let $G=K_{3,3,3}\vsum K_2$.
By Theorem~\ref{thm:cut-vertex},
\begin{align*}
  \mr(F,G)&=\min\{\mr(F,K_{3,3,3})+\mr(F,K_2),\,
  \mr(F,K_{3,3,2})+\mr(F,K_1)+2\}\\
  &=\min\{3+1,\,2+0+2\}=4.
\end{align*}
But since for either of the two nonisomorphic graphs $K_{3,3,2}\vsum
K_2$ arising from different choices of $v$, we have
$\mr(F,K_{3,3,2}\vsum K_2)\leq\mr(F,K_{3,3,2})+\mr(F,K_2)=2+1=3$ by
Corollary~\ref{cor:vsum-no-more-than-sum}, it follows that
$K_{3,3,3}\vsum K_2 \in \F_4(F)$.  This graph did not occur in the table
$\F_4(\Fld_2)$ because $\mr(\Fld_2,K_{3,3,3})=2$
\cite{barrett-vdHL-minrank2-infinite}.  It is also easy to see that
$K_{3,3,3}\cup K_2\in\F_4(F)$ if $\character F\neq 2$.  However, it is
difficult to analyze the structure of all graphs in $\F_4(F)$ with
vertex connectivity less than 2.  It is difficult to understand even
the graphs in $\F_4(F)$ that are of the form $G\vsum K_2$.  Sometimes
$G\in\F_3(F)$, but frequently it is not.  We do know, however, that
$\F_4(F)$ is infinite if $F$ is $\mathbb{R}$ or $\mathbb{C}$
\cite{hall-minrank}.

In examining the list of graphs in Section~\ref{sec:graphs-in-F_4}, we
see that some of the bounds obtained in Sections~\ref{sec:3k2-or-p3Vp3}--\ref{sec:p4-lemma} do not appear
to be sharp.  For instance, there is no graph in
Section~\ref{sec:graphs-in-F_4} with 8 vertices that has an induced $P_4$,
even though the bound in Lemma~\ref{lem:p4} is 8 vertices.  This is
because there are graphs in $\F_4(\Fld_2,P_4)$ that are not in
$\F_4(\Fld_2)$.
\begin{figure}
  \centering
  \input{graphs/p4nottight.tex}
  \caption{An 8 vertex graph in $\F_4(\Fld_2,P_4)\setminus \F_4(\Fld_2)$.}
\label{fig:p4-not-tight}
\end{figure}
For example, Figure~\ref{fig:p4-not-tight} shows a graph on 8
vertices which is in $\F_4(\Fld_2,P_4)$ (when the induced $P_4$
contains both center vertices), as can be checked by hand or by using
the Magma functions in the appendix.  However, the graph is not in
$\F_4(\Fld_2)$, since deleting one of the center vertices yields
graph~56 in Section~\ref{sec:graphs-in-F_4}.  This shows that
Lemma~\ref{lem:p4} does indeed provide a sharp bound for the number of
vertices in a graph in $\F_4(\Fld_2,P_4)$.

We have succeeded in obtaining a sharp bound on the number of vertices
in a minimal forbidden subgraph for the class of graphs whose minimum
rank is at most 3 over $\Fld_2$.  We have also generated a complete
list of these minimal forbidden subgraphs, thereby giving a structural
characterization for the graphs having minimum rank 4 or more over
$\Fld_2$.  Since this result leads to a method for generating or
recognizing all such graphs, it also leads to a theoretical procedure for
determining whether a given graph has minimum rank at most 3 over
$\Fld_2$.

\appendix

\section{Magma programs}
\label{sec:magma-programs}

\begin{verbatim}
// We are working in F_2.
F:=FiniteField(2);

// This function returns all matrices in S(F_2,G) by adding
// all possible diagonal matrices to the adjacency matrix of G.
matrices_in_S:=function(graph)
    return {DiagonalMatrix(F,x)+AdjacencyMatrix(graph): 
            x in Subsequences({x: x in F}, #Vertices(graph))};
end function;

// This function returns the minimum rank of a matrix by brute
// force computation.
minrank:=function(graph)
    return Min({Rank(m): m in matrices_in_S(graph)});
end function;

// This function returns the matrices in S(F_2,G) that attain
// the minimum rank.
minrank_matrices:=function(graph)
    return {m: m in matrices_in_S(graph) | Rank(m) eq minrank(graph)};
end function;

// This function returns true if and only if a subgraph of graph is
// isomorphic to a graph in graphlist 
// (i.e., if graph is forbidden by graphlist).
isomorphic_subgraph:=function(graph,graphlist)
    if exists(t){<subgraph,fgraph>:
            subgraph in {sub<graph|s>: s in Subsets(Set(VertexSet(graph)))},
            fgraph in graphlist 
            | IsIsomorphic(subgraph,fgraph)} then
         return true;
    else
        return false;
    end if;
end function;

// This is another version of the isomorphic_subgraph function.
isomorphic_subgraph:=function(graph,graphlist)
    for subgraph in {sub<graph|s>: s in Subsets(Set(VertexSet(graph)))} do
        if exists(t){ fgraph: fgraph in graphlist | 
                      IsIsomorphic(subgraph,fgraph)} then
             return true;
        end if;
    end for;
    return false;
end function;

// This function appends a list of forbidden subgraphs with
// numvertices vertices to forbiddengraphs.  The geng program
// must be in the current directory.
generate_forbidden_graphs:=function(numvertices,forbiddengraphs)
    allgraphs:=OpenGraphFile("cmd geng "
        *IntegerToString(numvertices), 0, 0);
    while true do
        more, graph:=NextGraph(allgraphs);
        if more then
            if minrank(graph) ge 4 
            and not isomorphic_subgraph(graph,forbiddengraphs) then
                Include(~forbiddengraphs,graph);
            end if;
        else
            break;
        end if;
    end while;
    return forbiddengraphs;
end function;

// Initialize the forbiddengraphs set and generate the forbidden
// subgraphs with 8 or fewer vertices.
forbiddengraphs:={};
for i in [1..8] do
    forbiddengraphs:=generate_forbidden_graphs(i,forbiddengraphs);
end for;

// Now forbiddengraphs contains all graphs in \mathcal{F}_4(F_2) as
// Magma graphs.
\end{verbatim}

\section{Field independent proof of Theorem~\ref{thm:cut-vertex}}
\label{appendix:cut-vertex}

First recall a definition, a well-known fact, and the statement of the theorem.
\begin{definition*}
  Let $G$ and $H$ be graphs on at least two vertices, each having a
  vertex labeled $v$.  Then $G\vsum H$ is the graph obtained from $G\cup
  H$ by identifying the two vertices labeled $v$.
\end{definition*}

\begin{lemma*}[{\cite{nylen-minrank}}]
  If $F$ is any field and $G$ is a graph with a vertex $v$, then
  $\mr(F,G-v)\leq\mr(F,G)\leq\mr(F,G-v)+2$.
\end{lemma*}

\begin{theorem*}[{\cite{hsieh-minrank,BFH1-minrankpath}}]
  Let $F$ be any field and let $G$ and $H$ be graphs on at least two
  vertices, each having a vertex labeled $v$.  Then
\begin{equation}
  \mr(F,G \vsum H)=\min\{\mr(F,G)+\mr(F,H),\, \mr(F,G-v)+\mr(F,H-v)+2\}.\label{eq:1}
\end{equation}
\end{theorem*}

\begin{proof}
  Since $v$ is a cut vertex of the connected graph $G\vsum H$,
  $(G\vsum H)-v=(G-v)\cup(H-v)$.  By the lemma and
  Observation~\ref{obs:other-fields}, 
  \begin{equation*}
    \mr(F,G\vsum H)\leq\mr(F,G-v)+\mr(F,H-v)+2.
  \end{equation*}

Let $v$ be the last vertex of $G$ and the first vertex of $H$.  Let 
\begin{equation*}
M=
\begin{bmatrix}
  A&b\\
  b^T&c_1
\end{bmatrix}\in S(F,G) \quad \text{ and } \quad 
N=
\begin{bmatrix}
  c_2&d^T\\
  d&E
\end{bmatrix}\in S(F,H),
\end{equation*}
such that $\rank M=\mr(F,G)$ and $\rank N=\mr(F,H)$.  Let
\begin{equation*}
  \hat M =
  \begin{bmatrix}
    A&b&0\\
    b^T&c_1&0\\
    0&0&0
  \end{bmatrix}\quad \text{ and } \quad
\hat N =
\begin{bmatrix}
  0&0&0\\
  0&c_2&d^T\\
  0&d&E
\end{bmatrix}.
\end{equation*}
Then $\hat M + \hat N \in S(F,G\vsum H)$ so
\begin{align*}
  \mr(F,G\vsum H) &\leq \rank(\hat M+\hat N)\\
  &\leq \rank \hat M + \rank \hat N
  = \rank M + \rank N \\
  &= \mr(F,G)+\mr(F,H).
\end{align*}

This proves the $\leq$ in \eqref{eq:1}.
\bigskip

Now let $M\in S(F,G\vsum H)$ with $\rank M=\mr(F,G\vsum H)$.  Write
\begin{equation*}
  M=
  \begin{bmatrix}
    A&b&0\\
    b^T&c&d^T\\
    0&d&E
  \end{bmatrix}.
\end{equation*}

Now
\begin{align}
  \label{eq:2} \rank A + \rank E &\leq \rank \bmat A&b&0\\0&d&E \emat \\
  \label{eq:3}&\leq \rank M\\
  \label{eq:4}&\leq \rank A + \rank E + 2.
\end{align}
It follows that one of the three inequalities \eqref{eq:2}, \eqref{eq:3},
or \eqref{eq:4} is an equality.

\begin{enumerate}[{\listindent}I.]
\item \label{item:3}Suppose that \eqref{eq:2} and \eqref{eq:4} are strict
  inequalities.  Then
  \begin{equation*}
    \rank M = \rank \bmat A&b&0\\0&d&E \emat = \rank A + \rank E + 1.
  \end{equation*}
  Consequently $\bmat b\\d \emat \not\in \col\bmat A&0\\0&E \emat$, so
  either $b\not\in\col(A)$ or $d\not\in\col(E)$.  Assume $b\not\in
  \col(A)$. Then $b^T\not\in \row(A)$, so
  \begin{equation*}
    \rank M = \rank \bmat A&b&0\\b^T&c&d^T\\0&d&E \emat > \rank \bmat
    A&b&0\\0&d&E \emat, 
  \end{equation*}
  a contradiction.  Therefore, this case does not occur.  So either
  \eqref{eq:2} or \eqref{eq:4} is an equality.

\item\label{item:4} Suppose \eqref{eq:2} is an equality.  Then
  \begin{equation*}
    \rank \bmat A&0\\0&E \emat = \rank \bmat A&0&b\\0&E&d \emat.
  \end{equation*}
  Thus $\bmat b\\d \emat \in \col \bmat A&0\\0&E \emat$, which implies
  that $b=Au,\, d=Ev$ for some vectors $u$ and $v$.  Then
  \begin{equation*}
    \hat A=\bmat A&Au\\u^TA&u^TAu \emat = \bmat A&b\\b^T&u^TAu \emat \in S(F,G)
  \end{equation*}
  and $\rank \hat A=\rank A$.  Similarly,
  \begin{equation*}
    \hat E=\bmat v^TEv&v^TE\\Ev&E \emat = \bmat v^TEv&d^T\\d&E \emat \in S(F,H)
  \end{equation*}
  and $\rank \hat E=\rank E$.  It follows that
  \begin{align*}
    \mr(F,G\vsum H)&=\rank M\\
    &\geq\rank \bmat A&0\\0&E \emat = \rank A+\rank E = \rank \hat A +
    \rank \hat E\\
    &\geq \mr(F,G)+\mr(F,H).
  \end{align*}
\item\label{item:5} Suppose that \eqref{eq:4} is an equality.  Since $A\in S(F,G-v)$
  and $B\in S(F,H-v)$, $\rank A\geq \mr(F,G-v)$ and $\rank B \geq \mr(F,H-v)$.
  Then
  \begin{equation*}
    \mr(F,G\vsum H)=\rank M \geq \mr(F,G-v)+\mr(F,H-v)+2.
  \end{equation*}
\end{enumerate}

Combining cases \ref{item:3}, \ref{item:4}, and \ref{item:5}, we have
proven the $\geq$ in \eqref{eq:1}.
\end{proof}

\section{Sage code to generate forbidden graphs}
This appendix contains a translation of the code in
Appendix~\ref{sec:magma-programs} for Sage (see
\url{http://www.sagemath.org}).

\begin{lstlisting}
# This code is written for Sage 3.1.1
# See http://www.sagemath.org

def matrices_in_S(graph):
    """
    Return all matrices in S($\mathbb{F}_2$,graph) by adding all
    possible diagonal matrices to the adjacency matrix of the graph.
    The matrices are returned as matrices over $\mathbb{F}_2$.
    """
    F = FiniteField(2)
    A = graph.adjacency_matrix().change_ring(F)
    return [A+diagonal_matrix(list(diagonal)) 
            for diagonal in VectorSpace(F,graph.order())]


def minrank(graph):
    """
    Return the minimum rank of a graph over $\mathbb{F}_2$ by
    exhaustive enumeration.
    """
    return min(m.rank() for m in matrices_in_S(graph))

def minrank_matrices(graph):
    """
    Return the matrices in S($\mathbb{F}_2$,graph) that attain the
    minimum rank of G.
    """
    minimumrank = minrank(graph)
    return [m for m in matrices_in_S(graph) if m.rank() == minimumrank]

def isomorphic_subgraph(graph, graphlist):
    """
    Return True if a subgraph of graph is isomorphic to a graph in
    graphlist, otherwise return False.
    """
    # Get the numbers of vertices in graphlist
    n = graph.order()
    subgraph_sizes = set([g.order() for g in graphlist if g.order() <= n])
    for i in subgraph_sizes:
        graphlist_check = [g for g in graphlist if g.order() == i]
        for vertices in Combinations(graph.vertices(), i):
            subgraph = graph.subgraph(vertices)
            if any(subgraph.is_isomorphic(g) for g in graphlist_check):
                return True
    return False


def generate_forbidden_graphs(numvertices, forbiddengraphs):
    """
    Return all graphs having numvertices vertices and minimum rank at
    least 4 that do not contain a subgraph isomorphic to a graph in
    forbiddengraphs.
    """
    return [g for g in graphs(numvertices) 
            if minrank(g) >= 4 and not isomorphic_subgraph(g, forbiddengraphs)]


# Wait for a while for these next few commands to complete.
forbiddengraphs=[]
for i in [1..8]:
    forbiddengraphs += generate_forbidden_graphs(i, forbiddengraphs)

# Now forbiddengraphs contains all graphs in $\mathcal{F}_4(\mathbb{F}_2)$ as
# Sage graphs.
\end{lstlisting}

\bibliographystyle{alpha}
\bibliography{}

\end{document}